\documentclass[a4paper,11pt,reqno,twosided]{amsart} 
\usepackage{amsmath,amsthm,amssymb,upref}
\usepackage[mathscr]{eucal}
\usepackage{graphpap,latexsym,epsf}
\usepackage{epsf}
\usepackage{color,graphics}
\usepackage[numbers,sort&compress,longnamesfirst]{natbib}

\theoremstyle{plain}
\newtheorem{theorem}{Theorem}[section]  
\newtheorem{corollary}[theorem]{Corollary}  
\newtheorem{lemma}[theorem]{Lemma}  

\newtheorem*{theorem*}{Theorem}

\theoremstyle{plain}  
\newtheorem{remark}[theorem]{Remark}

\newtheoremstyle{citing}
  {3pt}
  {3pt}
  {\itshape}
  {}
  {\bfseries}
  {.}
  {.5em}
  {\thmnote{#3}}

\theoremstyle{citing}

\numberwithin{equation}{section}

\newlength{\intwidth}
\makeatletter
\DeclareRobustCommand{\cpvint}[2]
    {\mathop{%
       \text{%
         \settowidth{\intwidth}{%
           \ifx\ilimits@\displaylimits
             $\int_{#1}^{#2}$%
           \else
             $\int$%
           \fi}%
         \makebox[0pt][l]{\makebox[\intwidth]{$\text{C}$}}%
         $\int_{#1}^{#2}$}}}
\makeatother

\makeatletter
\DeclareRobustCommand{\cpvintsmall}[2]
    {\mathop{%
       \text{%
         \settowidth{\intwidth}{%
           \ifx\ilimits@\displaylimits
             $\int_{#1}^{#2}$%
           \else
             $\int$%
           \fi}%
         \makebox[0pt][l]{\makebox[\intwidth]{$\text{{\tiny C}}$}}%
         $\int_{#1}^{#2}$}}}
\makeatother

\renewcommand{\a }{\alpha }

\newcommand{\Di }{\mathcal{D}^{1,2}(\R^N) } 
\newcommand{\Did }{\mathcal{D}^{1,2}(\R^3) } 
\newcommand{\stereo }{\mathcal{S} }

\renewcommand{\l }{\lambda }

\newcommand{\R}{\mathbb{R}}

\newcommand{\charak}{\chi}
\newcommand{\dist}{\text{dist}}

\newcommand{\xpfeil}{\xrightarrow}  
  
\newcommand{\rand}{\partial} 
\newcommand{\where}{\,:\:}

\newcommand{\sgn}{\text{sgn}}
  
\newcommand{\supm}{\mathop{\sup}\limits}
  
\newcommand{\cupl}{\mathop{\cup}\limits}

\newcommand{\intg}{\mathop{\int}\limits}  
\newcommand{\sumg}{\mathop{\sum}\limits}

\newcommand{\laplace}{\Delta}

\newcommand{\di}{\;d}  
  
\newcommand{\nz}{{\mathbb N}}

\newcommand{\rz}{{\mathbb R}}

\renewcommand{\phi}{\varphi}

\begin{document}
 
\title[] {A priori estimates for the scalar curvature equation on $S^3$}
\author{Matthias Schneider}
\address{Ruprecht-Karls-Universit\"at \\
         Mathematisches Institut \\ 
         Im Neuenheimer Feld 288 \\
         69120 Heidelberg, Germany}
\email{mschneid@mathi.uni-heidelberg.de} 
\date{December 1, 2004}  
\keywords{prescribed scalar curvature, a priori estimates}
\subjclass{35J60, 35J20, 55C21}

\begin{abstract}
We obtain a priori estimates for solutions to the prescribed scalar curvature
equation on $S^3$. The usual non-degeneracy assumption on the curvature function is
replaced by a new condition, which is necessary
and sufficient for the existence of a priori estimates, when the curvature
function is a positive Morse function. 
\end{abstract}

\maketitle

\section{Introduction}
Let $N\ge 3$ and $S^N$ be the standard sphere with round metric $g_0$ induced by $S^N = \rand B_1(0) \subset \rz^{N+1}$.
We study the problem: 
Which functions $K$ on $S^N$ occur as scalar curvature of metrics 
$g$ conformally equivalent to $g_0$? 
Writing $g= \varphi^{4/(N-2)} g_0$ this is equivalent to solving (see \cite{Aubin98})
\begin{align}
\label{eq:30}
-\frac{4(N-1)}{N-2}\laplace_{S^N}\varphi + N(N-1)\varphi= K(\theta) (\varphi)^{\frac{N+2}{N-2}},\;  \varphi>0  
\ \text{ in } S^N.
\end{align} 
In stereographic coordinates $\stereo_{\theta}(\cdot)$ centered at some point $\theta \in S^N$ equation \eqref{eq:30}
is equivalent to
\begin{align}
\label{eq:47}
-\laplace u = \frac{K\circ\stereo_{\theta} (x)}{N(N-1)} u^{\frac{N+2}{N-2}},\; u>0  \ \text{ in } \rz^N,
\end{align}
where
\begin{align}
\label{eq:116}
u(x) = {\mathcal R}_\theta(\phi)(x):= (N(N-2))^{\frac{N-2}{4}}(1+|x|^2)^{-\frac{N-2}{2}} \phi \circ \stereo_{\theta}(x).  
\end{align}
Obviously, to solve \eqref{eq:30} the function $K$ has to be positive somewhere. Moreover, there 
are the {\em Kazdan-Warner} obstructions \cite{KazdanWarner75,BourguignonEzin87}, if $\phi$ solves \eqref{eq:30} then
\begin{align*}
\int_{S^N} \nabla x_j \cdot \nabla K\, \varphi^{\frac{2N}{N-2}} =0 \ \text{ for }j=1 \dots N+1.  
\end{align*}
In particular, a monotone function of $x_1$ can not be realized as the scalar curvature of a metric conformal to
$g_0$.\\ 
Numerous studies have been made on equation \eqref{eq:30} and various sufficient conditions
for its solvability have been found (see \cite{AmbAzPer99,YYLi96,YYLi95,ChenLin01,ChenLi01,Bianchi96,AubinBahri97}
and the reference therein), usually under a non-degeneracy assumption on $K$. On $S^3$ a positive function
$K$ is non-degenerate, if
\begin{align}
\label{eq:118}
\tag{nd}
\laplace_{S^3} K(\theta) \neq 0 \text{ if } \nabla K(\theta) =0.   
\end{align}   
For positive Morse functions $K$ on $S^3$ it is shown in 
\cite{SchoenZhang96,BahriCoron91,ChangGurskyYang93} 
that \eqref{eq:30} is solvable if $K$ satisfies \eqref{eq:118} and
\begin{align}
\label{eq:34}
d:= -\Big( 1+\sum_{\substack{\nabla K(\theta)=0,\\ \laplace_{S^3} K(\theta)<0}} (-1)^{\text{ind}(\theta)}\Big) \neq 0,
\end{align}
where $\text{ind}(\theta)$ is the Morse index of $K$ at $\theta$. We are interested in the case
when $N=3$ and the non-degeneracy assumption \eqref{eq:118} is not satisfied.\\     
To obtain the existence result \citet{BahriCoron91} use a detailed analysis of the gradient flow of \eqref{eq:30}
and \citet{SchoenZhang96} approximate \eqref{eq:30} by subcritical problems $p \nearrow \frac{N+2}{N-2}$, 
which are always solvable, and analyze the possible blow-up of solutions.\\ 
We follow the approach
suggested in \cite{ChangGurskyYang93} and use a continuity method. We join 
the curvature function $K$ to the constant function $K_0\equiv 6$
by a one parameter family $K_t(\theta) := 6(1+t k(\theta))$, where 
$
k(\theta):= \frac{1}{6}(K(\theta)-6),  
$
and consider
\begin{align}
\label{eq:eq1}
 -8\laplace_{S^3}\varphi + 6\varphi= 6(1+t k(\theta)) \varphi^{5},\;  \varphi>0 \ \text{ in } S^3, 
\end{align}
or in stereographic coordinates using \eqref{eq:116} and $k_\theta(x) := k \circ \stereo_\theta(x)$
\begin{align}
\label{eq:eq2}
-\laplace u = (1+t k_\theta(x)) u^{5} \text{ in }\rz^3, \; u>0. 
\end{align}
In general there are no a priori $L^\infty$-estimates for \eqref{eq:eq1} or \eqref{eq:30}
due to the noncompact group of conformal transformations
of $S^N$ acting on solutions: the solutions of \eqref{eq:47} for $k\equiv 0$ form a noncompact manifold 
(see \cite{CaffarelliGidasSpruck89,GidasNiNirenberg81})
\begin{align*}
Z:= \big\{z_{\mu,y}(x):= \mu^{-\frac{N-2}{2}}(N(N-2))^{\frac{N-2}{4}} 
&\Big(1+\big|\frac{x-y}{\mu}\big|^2\Big)^{-\frac{N-2}{2}}\\ 
&\where y \in \rz^N, \mu>0\big\},    
\end{align*}
where $z_{\mu,y}(y) \to \infty$ as $\mu \to 0$.\\
\citet*{ChangGurskyYang93} show that if $K\in C^2(S^3)$ is positive and satisfies \eqref{eq:118}
then for every $\delta>0$ there is a constant $C=C(\delta,K)>0$
such that for all $t\in [\delta,1]$ and solutions $\phi_t$ of \eqref{eq:eq1} we have
\begin{align*}
C^{-1} \le \phi_t(\theta) \le C  \ \text{ and } \|\phi_t\|_{C^{2,\a}(S^3)}\le C.
\end{align*}
Furthermore, they compute the Leray-Schauder degree for \eqref{eq:eq1} for $t>0$ small, and show 
that it equals $d$ in \eqref{eq:34} if $K$ is a Morse function. The a priori estimate implies 
the invariance of the degree as the parameter $t$ moves to $1$ and gives a solution to \eqref{eq:eq1}
if $d\neq0$. \citet{ChenLin01} show that if $K\in C^2(S^3)$ is a non-degenerate Morse function then $C$ 
may be chosen independently of $\delta>0$.\\
Hence, if \eqref{eq:118} fails, we face two problems: Is the a priori bound still valid and 
how do critical points of $K$ with $\laplace_{S^3} K =0$ occur in the index count condition \eqref{eq:34}.
Here, we will mainly deal with the question about the a priori bound of solutions.\\  
In the following, unless otherwise stated, we will always assume $N=3$ and that the function 
$K \in C^5(S^3)$ is positive.
To give our main results we need the following notation.
For $k \in C^5(S^3)$ we write $k_\theta= k \circ \stereo_\theta$ and
for a critical point $\theta$ of $k$ we let 
\begin{align}
\label{eq:31}
\begin{split}
a_0(\theta)&:= \cpvint {\rz^3}{} \Big(k_{\theta}(x)-\sum_{\ell =0}^2 \frac{1}{\ell !} D^\ell k_{\theta}(0)(x)^\ell\Big)|x|^{-6}, \\
a_1(\theta)&:=  \laplace^2k_{\theta}(0)
+\nabla(\laplace k_{\theta}(0)) \cdot \big(D^2k_{\theta}(0)\big)^{-1}\nabla(\laplace k_{\theta}(0)),\\
a_2(\theta)&:= k_{\theta}(0)a_1(\theta)
- \frac{15}{8\pi}\int_{\rand B_1(0)}\big|D^2k_{\theta}(0)(x)^2\big|^2,  
\end{split}
\end{align}  
where all differentiations are done in $\rz^3$ and $\cpvintsmall{}{}$ is the Cauchy principal value of the integral,
\begin{align*}
\cpvint{\rz^3}{} f(x) := \lim_{r \to 0} \int_{\rz^3\setminus B_r(0)}f(x). 
\end{align*}
\begin{theorem}
\label{t:apriori}
Suppose $1+k\in C^5(S^3)$ is positive and satisfies 
\begin{align*}
D^2 k_\theta(0) \text{ is invertible, if }  
\theta \in \mathcal{A} := \{\theta \in S^3 \where \nabla k(\theta)=0 \text{ and }\laplace_{S^3} k(\theta)=0\}. 
\end{align*}
Thus, $\mathcal{A}$ is discrete and hence finite. Denote by $M$ the finite set
\begin{align*}
M:= \big\{\theta \in S^3 \where \theta \in \mathcal{A}, \, a_0(\theta)= 0, \,\text{and } a_2(\theta)\neq 0\big\}. 
\end{align*}
Then for every $\delta>0$ there is a constant $C=C(k,\delta)>0$ such that for all 
\begin{align*}
t \in (0,1]\setminus \cupl_{\theta \in M} B_\delta\Big(-\frac{a_1(\theta)}{a_2(\theta)}\Big)  
\end{align*}
and solutions $\phi_t$ of \eqref{eq:eq1} we have
\begin{align*}
C^{-1} \le \phi_t(x) \le C  \ \text{ and } \|\phi_t(x)\|_{C^{2,\a}(S^3)}\le C.
\end{align*}
\end{theorem}
Theorem~\ref{t:apriori} extends the known a priori estimates to the case when
\eqref{eq:118} may fail. If $k$ satisfies \eqref{eq:118} the solutions are uniformly bounded with respect to $t \in(0,1]$.
Moreover, we get uniform estimates for $t \in(0,1]$ if 
$M^*=\emptyset$, where
\begin{align*}
M^*:= \{\theta \in M \where 0\le-a_1(\theta)/a_2(\theta)\le 1\}.  
\end{align*}
Our results are optimal since we construct for every $\theta \in M^*$ solutions $\phi_t$ which blow up
as $t\to -a_1(\theta)/a_2(\theta)$. We say that $(t_i,\phi_i)$ {\em blow up } at the 
{\em blow-up point} $\theta \in S^3$,
if $\phi_i$ solves \eqref{eq:eq1} with $t=t_i$, the sequence $(t_i)$ is bounded, 
and there is $(\theta_i)$ converging to $\theta$
such that $\phi_i(\theta_i) \to \infty$.
\begin{theorem}
\label{t:a priori_existence_plus}
Under the assumptions of Theorem \ref{t:apriori}  let 
$$M^*_+ := \{\theta \in M \where 0<-a_1(\theta)/a_2(\theta)\le 1\}.$$
Then there is $\delta>0$ such that for any $\theta \in M^*_+$ there exists a unique $C^1$-curve
\begin{align*}
\{0<\mu<\delta\} \ni \mu \mapsto (t^\theta(\mu),\phi^\theta(\mu,\cdot)) 
\in (\delta,1+\delta) \times C^{2,\a}(S^3), 
\end{align*}
such that as $\mu \to 0$
\begin{align*}
t^\theta(\mu) =-\frac{a_1(\theta)}{a_2(\theta)}+ O(\mu^{\frac14}),  
\end{align*}
and $\phi^\theta(\mu,\cdot)$ solves \eqref{eq:eq1} for $t=t^\theta(\mu)$ and blows up like
\begin{align*}
\|{\mathcal R}_\theta(\phi^\theta(\mu,x))-(1+t^\theta(\mu)k(\theta))^{-\frac14} 
z_{\mu,0}(x)\|_{\Did\cap C^2(B_1(0))} = O(\mu^2).  
\end{align*} 
The curves are unique, in the sense that, 
if $(t_i,\phi_i) \in (\delta,1+\delta) \times C^{2,\a}(S^3)$ blow
up at some $\theta \in S^3$ then $\theta \in M^*_+$ and there is a sequence of positive numbers 
$(\mu_i)$ converging to zero such that 
$(t_i,\phi_i) = (t^\theta(\mu_i),\phi^\theta(\mu_i,\cdot))$
for all but finitely many $i \in \nz$. 
\end{theorem}
Hence, for Morse functions we obtain 
\begin{corollary}
\label{c.a priori_morse}
Suppose $1+k\in C^5(S^3)$ is a positive Morse function. There exists $\delta_0>0$,
such that for any $0<\delta<\delta_0$ the solutions of \eqref{eq:eq1} are uniformly bounded for $t\in [\delta,1+\delta]$,
if and only if $M^*_+=\emptyset$.
\end{corollary}
For $\theta \in M$ with $a_1(\theta)=0$ there is always the trivial curve of solutions,
\begin{align*}
\mu \mapsto \big(0,(\mathcal R_\theta)^{-1} z_{\mu,0}\big) \in \rz \times C^{2,\a}(S^3),  
\end{align*}
which blow up at $\theta$ as $\mu \to 0$. In order to find a nontrivial curve, i.e. $t(\mu) \in \rz\setminus\{0\}$,
we need to consider
\begin{align}
\label{eq:33}
\begin{split}
a_3(\theta)&:= 
\frac{12}{\pi^2} \Big(D^2k_\theta(0)\Big)^{-1} \nabla(\laplace k_\theta(0)) \cdot 
\cpvint{\rz^3}{} \Big(\nabla k_\theta(x)-  T_{\nabla k_\theta,0}^2(x)\Big)|x|^{-6}\\
&\quad +\frac{48}{\pi^2} \big(D^2k_\theta(0)\big)^{-1} \nabla(\laplace k_\theta(0))
\cdot
\cpvint{\rz^3}{} \Big(k_\theta(x)-T_{k_\theta,0}^3(x)\Big)\frac{x_i}{|x|^{8}}\\
&\quad -\frac{120}{\pi^2} 
\cpvint{\rz^3}{} \Big(k_\theta(x)- T_{k_\theta,0}^4(x)\Big) \frac1{|x|^8},
\end{split}
\end{align}
where we abbreviate the $m$th Taylor polynomial of $k$ in $y$ by
\begin{align*}
T_{k,y}^m(x) := \sum_{\ell=0}^m\frac{1}{\ell!}D^\ell k(y)(x-y)^\ell.
\end{align*}
\begin{theorem}
\label{t:a priori_existence_null}
Under the assumptions of Theorem \ref{t:apriori} suppose $k\in C^6(S^3)$ and let
\begin{align*}
M^*_0 := \{\theta \in M \where a_1(\theta)=0 \text{ and } a_3(\theta)\neq 0\}.   
\end{align*}
Then there is $\delta>0$ such that
for any  $\theta \in M^*_0$ there exists a unique $C^1$-curve
\begin{align*}
\{0<\mu<\delta\} \ni \mu \mapsto (t(\mu),\phi(\mu,\cdot)) 
\in \big((-\delta,1+\delta)\setminus\{0\}\big) \times C^{2,\a}(S^3), 
\end{align*}
such that as $\mu \to 0$
\begin{align*}
t(\mu) =-\frac{a_3(\theta)}{a_2(\theta)}\mu + O(\mu^{1+\frac14}), 
\end{align*}
and $\phi(\mu,\cdot)$ solves \eqref{eq:eq1} for $t=t(\mu)$ and blows up like
\begin{align*}
\|\mathcal{R}_\theta(\phi(\mu,x))- z_{\mu,0}(x)\|_{\Did\cap C^2(B_1(0))} = O(\mu^2).  
\end{align*} 
The curve is unique, in the sense that, if $(t_i,\phi_i) \in 
\big((-\delta,1+\delta)\setminus\{0\}\big) \times C^{2,\a}(S^3)$ blow
up at $\theta \in M^*_0$ then there is a sequence of positive numbers 
$(\mu_i)$ converging to zero such that 
$(t_i,\phi_i) = (t(\mu_i),\phi(\mu_i,\cdot))$
for all but finitely many $i \in \nz$. 
\end{theorem}
To illustrate our results we give an example. Suppose $k_\theta$ is given by 
\begin{align*}
k_\theta(x) =  1 + \frac{3x_1^2-2x_2^2-x_1^2}{(1+|x|^2)^2}
+ \frac{|x|^4}{(1+|x|^2)^3}\Big(b-\frac{a}{(1+|x|^2)}-\frac{1-a}{(1+|x|^2)^2}\Big).   
\end{align*}
Then $1+tk_\theta(x)$ is strictly positive for all $t\ge 0$, if $b \ge 0$ and $a\le 3$,  
and  $\laplace k_\theta(0) =0$ and $\nabla k_\theta(0)=\nabla \laplace k_\theta(0) =0$.
Furthermore,
\begin{align*}
a_0(\theta) = -\frac{\pi^2}{64}(35-48b+5a),\quad 
&a_1(\theta) = 120 (b-1),\\
a_2(\theta) = 120 (b-1)-56,\quad 
&a_3(\theta) = 75\Big(6b+\frac78a-\frac{63}{8}\Big).  
\end{align*}
Our results show: $\theta$ is not a blow-up point, 
if $a_0(\theta)\neq 0$, that is $a \neq -7+\frac{48}{5}b$, or $a_2(\theta)=0$, that is $b=\frac{22}{15}$.
Moreover, if $a_0(\theta)=0$ and $a_2(\theta)\neq 0$ then there is a 
curve of solutions $(t(\mu),\varphi(\mu,\cdot))$ which blow up at $\theta$ such that
\begin{align*}
t(\mu) = \frac{b-1}{\frac{22}{15}-b} + O(\mu^{\frac14}), \quad &\text{if } b \neq 1\\
t(\mu) = \frac{30}{56}\mu+ O(\mu^{1+\frac14}), \quad &\text{if } b=1.
\end{align*}

We sketch the strategy of the proofs of our main results and outline the remaining part of the paper.
The transformation in \eqref{eq:116} gives rise to a Hilbert space isomorphism between
$H^{1,2}(S^N)$ and $\Di$, 
where $\Di$ denotes the closure of $C_c^\infty(\rz^N)$ with respect to 
\begin{align*}
\|u\|^2:= \int_{\rz^N} |\nabla u|^2 =\langle u,u\rangle.   
\end{align*}
Due to elliptic regularity (see \cite{BrezisKato79}) and
Harnack's inequality it is enough to find a weak nonnegative solution of \eqref{eq:30} in $H^{1,2}(S^N)$, or of
the equivalent equation. Although we take advantage of both formulations, we mainly  
consider \eqref{eq:47}.
We use a finite dimensional reduction of Melnikov type developed in \cite{AmBa2,AmbAzPer99} and 
find solutions of \eqref{eq:47} as critical points of $f_t:\Di \to \rz$, where
\begin{align*}
f_t(u) := \frac{1}{2}\int_{\rz^N} |\nabla u|^2 - \frac{N-2}{2N} \int_{\rz^N} (1+t k(x))|u|^{\frac{2N}{N-2}}.
\end{align*}
For $t=0$ the functional $f_0$ possesses, as seen above, a $N+1$ dimensional
manifold of critical points $Z$. To setup the finite dimensional reduction we need to analyze
$Z$ and the spectrum of $f_0''(z)$ in detail, which is done for all $N\ge 3$ in Section \ref{s:preliminaries}.
For the rest of the paper we will only deal with the case $N=3$.
In Section \ref{s:blow_up} we recall without proof that if $N=3$ a sequence of solutions to \eqref{eq:eq1}
can only blow-up in a single point (see \cite{SchoenZhang96,YYLi95}) and fit this result into our framework.
Section \ref{s:expansion} contains the finite dimensional reduction of our problem. In contrast to 
\cite{AmbAzPer99}, where the reduction is performed for small $t$, we show that a finite dimensional
reduction of \eqref{eq:eq2} for large $t$ is still possible. 
We end up with a function $\vec{\a}: U \to \rz^4$, where $U \subset \rz \times Z$, such that the
zeros of $\vec{\a}(t,\cdot)$ correspond to solutions of \eqref{eq:eq2} with large $L^\infty$ norm.
We recall that $Z$ is parametrized by $\mu$ and $y$. 
Now, to construct or to rule out blow-up sequences it is enough to construct or exclude zeros
of $\vec{\a}(t,\cdot)$ for small $\mu$. To this end we need to expand $\vec{\a}$ up to order $5$ in $\mu$ and to 
compute derivatives of $\vec{\a}$, which is done in Sections \ref{s:expansion} and 
\ref{s:derivatives_alpha}. We see that $\theta$ can only be a blow-up point if $\nabla k(\theta)=0$ and
$\laplace k_\theta(0) =0$. In Section \ref{s:solv_alpha_equal_null} we finally obtain under the assumptions
of Theorems \ref{t:apriori}-~\ref{t:a priori_existence_null} that
there are $(t_i,\phi_i)$ which blow up at $\theta$ if and only if $\nabla k(\theta)=0$, $\laplace k_\theta(0)=0$, 
and there exist positive $(\mu_i)$ converging to $0$ such that
\begin{align*}
0= a_0(\theta) +\mu_i (a_1(\theta)+t_i a_2(\theta))+\mu_i^2 a_3(\theta)+ O(t_i\mu_i^{1+\frac14}+\mu_i^{2+\frac14}). 
\end{align*}
This gives our main results, which are stated and proved in Section \ref{s:a priori_estimates}.\\ 
In a subsequent paper \cite{schneider04} we use the above a priori estimates and compute 
under the assumptions of Theorem \ref{t:apriori}
the Leray-Schauder degree $d$ of the problem \eqref{eq:eq1}. We show that if $M^*=\emptyset$ and $k$
is a Morse function then,
\begin{align*}
d= -\Big(1-\sum_{\theta \in \text{Crit}(k)^{-}} (-1)^{\text{ind}(\theta)}\Big),\ \text{ where}
\end{align*}
\begin{align*}
\text{Crit}(k)^{-} := \big\{ \theta \in S^3 \where &\nabla k(\theta)=0 \text{ and }\\ 
&\lim_{\mu \to 0^+} \sgn\big(\laplace k_\theta(0)+a_0(\theta)\mu +a_1(\theta)\mu^2\big)= -1 \big\},    
\end{align*}
generalizing the existence result in \cite{SchoenZhang96,BahriCoron91,ChangGurskyYang93}.\\
Our approach yields information about blow-up sequences as precise as we want, that is of any order
in $\mu$ or $y$. For instance it is possible to compute the term of order $\mu^6$, which is
of interest when $a_3(\theta)$ is zero. But the necessary computations and terms,
as may already be seen in the expansion of order $5$, are getting rather bulky.
In higher dimensions $N\ge 4$ solutions may blow up in more than one point and our method, which still applies
with minor changes to $N\ge 4$, will only give information about ``one bubble'' blow-up.

\section{Preliminaries}
\label{s:preliminaries}
We define for $\mu>0$ and $y \in \rz^N$ the 
maps ${\mathcal U}_{\mu},\,{\mathcal T}_y: \Di \to \Di$ by
\begin{align*}
{\mathcal U}_{\mu}(u):= \mu^{-\frac{N-2}{2}} u\Big(\frac{\cdot}{\mu}\Big) \text{ and }
{\mathcal T}_{y}(u):= u(\cdot-y).  
\end{align*}
With this notation the critical manifold $Z$ is given by
\begin{align*}
Z = \{ z_{\mu,y}= {\mathcal T}_y \circ {\mathcal U}_\mu (z_{1,0}) \where y \in \rz^N,\, \mu>0\}.
\end{align*}
It is easy to check that the dilation ${\mathcal U}_{\mu}$ and the translation ${\mathcal T}_{y}$ conserve 
the norms $\|\cdot\|$ and the $L^{2^*}$-Norm $\|\cdot\|_{2^*}$, where $2^*:= 2N/(N-2)$. 
Thus for every $\mu>0$ and $y \in \rz^N$ 
\begin{align}
\label{eq:7}
\begin{split}
({\mathcal U}_{\mu})^{-1} = ({\mathcal U}_{\mu})^{t} = U_{\mu^{-1}},\
({\mathcal T}_{y})^{-1} &= ({\mathcal T}_{y})^{t} = {\mathcal T}_{-y}, \text{ and }\\ 
f_0 = f_0 \circ {\mathcal U}_{\mu} = f_0 \circ {\mathcal T}_{y}
\end{split}  
\end{align}
where $(\cdot)^{t}$ denotes the adjoint. 
Twice differentiating the identities for $f_0$ in \eqref{eq:7} yields 
\begin{align}
\label{eq:1}
f_0''(v) = 
({\mathcal T}_{y}\circ {\mathcal U}_{\mu})^{-1} \circ 
f_0''({\mathcal T}_{y}\circ {\mathcal U}_{\mu}(v)) \circ ({\mathcal T}_{y}\circ {\mathcal U}_{\mu}) \quad
\forall v \in \Di. 
\end{align}
Moreover, we see that $U(\mu,y,z):= {\mathcal T}_{y}\circ{\mathcal U}_{\mu}(z)$ maps 
$(0,\infty)\times \rz^N \times Z$
into $Z$, hence 
\begin{align}
\label{eq:2}
\begin{split}
\frac{\rand U}{\rand z}(\mu,y,z) = 
{\mathcal T}_{y}\circ{\mathcal U}_{\mu}:\: T_zZ \to T_{{\mathcal T}_{y}\circ{\mathcal U}_{\mu}(z)}Z 
\text{ and }\\ 
{\mathcal T}_{y}\circ{\mathcal U}_{\mu}:\: (T_zZ)^{\perp} \to (T_{{\mathcal T}_{y}\circ{\mathcal U}_{\mu}(z)}Z)^{\perp}.   
\end{split}
\end{align}
The tangent space $T_{z_{\mu,y}}Z$ at a point $z_{\mu,y} \in Z$ is spanned by $N+1$ orthonormal functions 
$\dot{\xi}_{\mu,y}^i$,
\begin{align*}
T_{z_{\mu,y}}Z = \langle \dot{\xi}_{\mu,y}^i \where i=0 \dots N \rangle,   
\end{align*}
where ${\dot\xi^i}_{\mu,y}$ denotes for $i=0$ the normalized tangent vector $\frac{d}{d\mu} z_{\mu,y}$
and for $1\le i \le N$ the normalized tangent vector 
$\frac{d}{d y_i} z_{\mu,y}=-\frac{\rand}{\rand x_i} z_{\mu,y}$.
By \eqref{eq:2} we obtain $${\dot\xi^i}_{\mu,y}= {\mathcal T}_y \circ {\mathcal U}_\mu ({\dot\xi^i}_{1,0}).$$
An explicit calculation gives for $1 \le i \le N$  
\begin{align*}
({\dot\xi}_{1,0})_i = 
\tau_{1} (1+|x|^2)^{-\frac{N}{2}}x_i,\,
\tau_{1}^2 := \frac{\Gamma(N+2)}{\pi^{N/2}\Gamma(2+{N}/{2})N}.    
\end{align*}
For $i=0$ we find
\begin{align*}
({\dot\xi}_{1,0})_0 &= \tau_{0}
(1+|x|^2)^{-\frac{N-2}{2}}\Big(1-\frac{2}{(1+|x|^2)}\Big),\,
\tau_{0}^2 := \frac{(N+1)\Gamma(N)}{\pi^{N/2}\Gamma({N}/{2})N (N+2)}. 
\end{align*}
Using the canonical identification of the Hilbert space $\Di$ with its
dual induced by the scalar-product 
we shall consider $f_t'(u)$ as an element of $\Di$ and $f_t''(u)$ as one of
$\mathcal{L}(\Di)$. With this identification $f_t''(u)$ is of the form $identity-compact$ 
(see \cite{AmbAzPer99}) and hence 
a Fredholm operator of index zero.\\ 
Since $f_0''(z_{\mu,y})$ is a self-adjoint, 
compact perturbation of the identity map in $\Di$, its spectrum $\sigma(f_0''(z_{\mu,y}))$
consists of point-spectrum, possibly accumulating at $1$. We fix $\lambda \in \sigma(f_0''(z_{\mu,y}))$ and a 
corresponding eigenfunction $u$. Then $u$ solves
\begin{align}
\label{eq:linearized}
-\laplace u - \frac{N+2}{N-2}z_{\mu,y}^{2^*-2}u = \lambda (-\laplace u).
\end{align}
We expand $u$ in spherical harmonics with center $y$
\[
u(y+r\vartheta)=\sum_{i=0}^{\infty}\sum_{l=1}^{c_i}v_{i,l}(r){Y}_{i,l}(\vartheta),\quad 
r\in\R^+,\quad\vartheta\in {\mathbb S}^{N-1},
\]
where 
\[v_{i,l}(r)=\int_{{\mathbb S}^{N-1}}u(y+r\vartheta)Y_{i,l}(\vartheta)\,d\vartheta, \quad
c_i:= \binom{N-1+i}{N-1} -\binom{N-3+i}{N-1},\] 
and $\{Y_{i,l}\}$
denote a $L^2(S^{N-1})$-basis of (real valued) spherical harmonics
satisfying for all $i \in \nz_0$  and $1\le l \le c_i$
\begin{align*}
-\Delta_{{\mathbb S}^{N-1}}{Y}_{i,l}=i(N+i-2)Y_{i,l}.
\end{align*}
There is a freedom in choosing such a $L^2$-basis and because the cases $i=1,2$ will be of special 
interest in the sequel we fix the basis-vectors in these cases. 
We set for $i=1$ and $1 \le l \le N$
\begin{align}
\label{eq:36}
{Y}_{1,l}\Big(\frac{x}{|x|}\Big) := \Big(\frac{2 \pi^{N/2}}{N\Gamma(N/2)}\Big)^{-1/2} \frac{x_l}{|x|}.   
\end{align}
For $i=2$ we introduce a more convenient notation and write 
\begin{align}
\label{eq:37}
\begin{split}
{Y}_{2,(l_1,l_2)}\Big(\frac{x}{|x|}\Big) &:=  
\Big(\frac{\pi^{N/2}}{\Gamma(2+N/2)}\Big)^{-1/2} |x|^{-2}\\
&\quad\times
\begin{cases}
\sqrt{2}\,x_{l_1} x_{l_2} & \text{if } 1 \le l_1<l_2 \le N\\
\sqrt{\frac{l_1-1}{l_1}}\, \Big(x_{l_1}^2-\frac{1}{l_1-1}\sumg_{m=1}^{l_1-1} x_m^2\Big) 
& \text{if } 2 \le l_1=l_2 \le N.  
\end{cases}  
\end{split}
\end{align}
Since $u$ solves (\ref{eq:linearized}) the functions ${v}_{i,l}$ satisfy for $i \in \nz_0$  and $1\le l \le c_i$
\begin{align*}
-{v}_{i,l}{''}-\frac{N-1}{r}\, {v}_{i,l}{'}+\frac{i(N+i-2)}{r^{2}}
 \,{v}_{i,l}=\frac{N+2}{(N-2)(1-\lambda)}z_{\mu,0}^{2^*-2}{v}_{i,l}.
\end{align*}
Making the transformation
\begin{align*}
v(r)= r^{-\frac{N-2}{2}} \zeta(\ln \mu+\ln r),  
\end{align*}
we obtain the equation for $i \in \nz_0$  and $1\le l \le c_i$
\begin{align*}
-{\zeta}_{i,l}{''}-\frac{N(N+2)}{4(1-\l)}\cosh^{-2}(t){\zeta}_{i,l}
=\left(-\bigg(\frac{N-2}{2}\bigg)^2-i(N+i-2)\right){\zeta}_{i,l}.
\end{align*}
Using the results in \cite[p. 74]{landau} or \cite{olver} as in \cite{FS03} 
we find 
\begin{align*}
\lambda_{i,j}= 1-\frac{N(N+2)}{\left((N+1)+2(i+j-1)\right)^2-1}.  
\end{align*}
The corresponding eigenfunction is given by 
\begin{align*}
\psi_{i,j}(t):= (1-\tanh(t)^2)^{\frac{N-2+2i}{4}} {\mathcal P}_j^{(\sigma_i,\sigma_i)}(\tanh(t)),  
\end{align*}
where ${\mathcal P}_j^{(\sigma,\sigma)}$ denotes the {\em Jacobi} polynomial defined in \eqref{eq:A1}
and $\sigma_i$ is given by
\begin{align*}
\sigma_i := \frac{N-2}{2}+i.  
\end{align*}
Consequently, $\sigma(f_0''(z_{\mu,y}))= \{\lambda_{i,j} \where i,j \in \nz_0\}$ and the eigenspace of the
eigenvalue $\lambda_{i,j}$ has dimension $c_i$ and is spanned by, $(l=1 \dots c_i)$
\begin{align}
\label{eq:115}
\begin{split}
\Phi^{\mu,y}_{i,j,l}(x) := a_{i,j}\; {\mathcal T}_{y}\circ {\mathcal U}_{\mu}
\Big(|x|^i& (1+|x|^2)^{-\frac{N-2}{2}-i}\\
&{\mathcal P}_j^{(\sigma_i,\sigma_i)}\big(1-2(1+|x|^2)^{-1}\big) Y_{i,l}\big(\frac{x}{|x|}\big)\Big),  
\end{split}
\end{align}
where the $a_{i,j}$ are given by
\begin{align*}
a_{i,j}^2 := \frac{2(N-1+2(i+j))j! \Gamma(N-1+2i+j)}{(N-2+2(i+j))(N+2(i+j)) \Gamma(N/2+i+j)^2}  
\end{align*}
to assure that the $\Phi^{\mu,y}_{i,j,l}$ are orthonormal. 
Since $Z$ is a manifold of critical points of $f_0'$, the tangent space $T_z Z$ at a point $z \in Z$ is contained 
in the kernel $N(f_0''(z))$ of $f_0''(z)$. As $\lambda_{i,j}=0$ if and only if $i+j=1$, 
the dimension of $N(f_0''(z))$ is $N+1$, which implies that
\begin{align}
\label{eq:15}
T_z Z=  N(f_0''(z)) \quad \text{ for all } z \in Z. 
\end{align}
More precisely, we have
\begin{align}
\label{eq:11}
\Phi^{\mu,y}_{1,0,l} = {\dot{\xi}}_{\mu,y}^{l} \text{ for } l=1,\dots, N \ \text{ and }\ 
\Phi^{\mu,y}_{0,1,1} = {\dot{\xi}}_{\mu,y}^{0}. 
\end{align}
If \eqref{eq:15} holds the critical manifold $Z$ is called non-degenerate (see \cite{AmBa2})
and the self-adjoint Fredholm operator 
$f_0''(z)$ maps the space $\Di$ into
$T_{z}Z^\perp$ and is invertible in ${\mathcal{L}}(T_{z}Z^\perp)$. 
From \eqref{eq:1} and \eqref{eq:2}, we obtain in this case 
\begin{align}
\label{eq:3}
\|(f_0''(z_{1,0}))^{-1}\|_{{\mathcal{L}}(T_{z_{1,0}}Z^\perp)} = 
\|(f_0''(z))^{-1}\|_{{\mathcal{L}}(T_{z}Z^\perp)} \quad \forall z \in Z.  
\end{align}

\section{Blow up analysis}
\label{s:blow_up}
Let $(K_i) \in C^1(S^3)$ satisfy for some $A_0$
\begin{align}
\label{eq:16}
 A_0^{-1} \le K_i(x) \le A_0 \text{ and } \|\nabla K_i\|_\infty \le A_0.
\end{align}
We have the following result (see \cite{SchoenZhang96,YYLi95})
\begin{theorem}
\label{t:blow_up_cite}
Suppose $(K_i) \in C^1(S^3)$ satisfies \eqref{eq:16} and $(\varphi_i) \in C^2(S^3)$ solves \eqref{eq:30} with $N=3$.
Then after passing to a subsequence either 
$(\varphi_i)$ is uniformly bounded in $L^\infty(S^3)$ and hence in $C^{2,\alpha}(S^3)$ by elliptic regularity or 
$(\varphi_i)$ has precisely one isolated simple blow-up point $\theta$, i.e. 
there exists a sequence $(\theta_i)$ of maxima of $\varphi_i$ converging to some $\theta \in S^3$
and $C=C(A_0)$ such that $\varphi_i(\theta_i)\to +\infty$ and in geodesic normal coordinates about $\theta_i$ given by
$\exp_{\theta_i}(\cdot)$ 
\begin{align*}
i\varphi_i(\theta_i)^{-2} \to 0,\\ 
\bigg\|\frac{\varphi_i\Big(\exp_{\theta_i}\Big(\frac{x}{\varphi_i(\theta_i)^{2}}\Big)\Big) }{\varphi_i(\theta_i)} 
-\Big(1+\frac{K_i(\theta_i)}{24}|x|^2\Big)^{-\frac12}\bigg\|_{C^{2,\alpha}(\{|x|\le 3i\})} \le i^{-4},\\
\varphi_i(x) \le C \varphi_i(\theta_i)^{-1}\dist_{S^3}(x,\theta_i)^{-1} \ \text{ for } \dist_{S^3}(x,\theta_i)\ge  i\varphi_i(\theta_i)^{-2}.
\end{align*}
\end{theorem}
We need a slightly different version of this result.
\begin{corollary}
\label{c:blow_up}
Under the assumptions of Theorem \ref{t:blow_up_cite} the sequence $(\varphi_i) \in C^2(S^3)$ is, after
passing to a subsequence, either uniformly bounded in $C^{2,\alpha}$ or 
there exist $\theta \in S^3$ and sequences $(\mu_i)\in (0,\infty)$, $(y_i) \in \rz^3$ satisfying
\begin{align*}
\lim_{i \to \infty} \mu_i= 0,\quad  \lim_{i \to \infty} y_i= 0,       
\end{align*}
such that in stereographic coordinates $\stereo_\theta(\cdot)$
about $\theta$ the function $u_i$ defined by the transformation \eqref{eq:116} 
satisfies
\begin{align}
\label{eq:119}
u_i - 6^{\frac14} \big(K_i\circ \stereo_{\theta}(y_i)\big)^{-\frac14} z_{\mu_i,y_i} 
\text{ is orthogonal to } T_{z_{\mu_i,y_i}}Z,\\
\label{eq:120}
\|u_i - 6^{\frac14} \big(K_i\circ \stereo_{\theta}(y_i)\big)^{-\frac14} z_{\mu_i,y_i}\|_{\Did} = o_{A_0}(1).
\end{align}
\end{corollary}
To prove the corollary we first need the following lemma, which is an easy consequence of Theorem \ref{t:blow_up_cite}.
\begin{lemma}
\label{l:blow_up}
Under the assumptions of Theorem \ref{t:blow_up_cite} the sequence $(\varphi_i) \in C^2(S^3)$ is, after
passing to a subsequence, either uniformly bounded in $C^{2,\alpha}$ or 
there exists a sequence $(\theta_i)$ of maxima of $\varphi_i$ converging to some $\theta \in S^3$
and such that $\varphi_i(\theta_i)\to +\infty$ and in stereographic coordinates $\stereo_\theta(\cdot)$
using the transformation \eqref{eq:116}
\begin{align*}
\|u_i - 6^{\frac14} K_i(\theta_i)^{-\frac14} z_{\mu_i,y_i}(x)\|_{\Did} = o_{A_0}(1),
\end{align*}
where $\mu_i \to 0$ and $y_i \to 0$ are given by
\begin{align*}
y_i:= \stereo_\theta^{-1}(\theta_i),\quad \mu_i:= \big(6/K_i(\theta_i)\big)^\frac12 \varphi_i(\theta_i)^{-2}.   
\end{align*}
\end{lemma}
\begin{proof}[Proof of Corollary \ref{c:blow_up}]
From Lemma \ref{l:blow_up} we infer that there are $\theta \in S^3$, $(\tilde{\mu_i})$, and $(\tilde{y_i})$
such that
\begin{align*}
\mu_i + |y_i|+
\|u_i - 6^{\frac14} \Big(K_i\circ \stereo_{\theta}(\tilde{y_i})\Big)^{-\frac14} z_{\tilde{\mu_i},\tilde{y_i}}\|_{\Did}
= o_{A_0}(1).
\end{align*}
For $i$ fixed, consider
\begin{align*}
d_i:= \inf_{\mu,y} \|u_i- 6^{\frac14} \Big(K_i\circ \stereo_{\theta}(\tilde{y_i})\Big)^{-\frac14} z_{\mu,y}\|_{\Did}^2. 
\end{align*}
Clearly $d_i= o_{A_0}(1)$ and therefore $d_i$ is attained at $\mu_i,y_i$ and
\begin{align*}
 u_i- 6^{\frac14} \Big(K_i\circ \stereo_{\theta}(\tilde{y_i})\Big)^{-\frac14} z_{\mu_i,y_i} \text{ is orthogonal to }
T_{z_{\mu_i,y_i}}Z.
\end{align*}
Since $z_{\mu_i,y_i}$ is orthogonal to $T_{z_{\mu_i,y_i}}Z$ relation \eqref{eq:119} follows. To prove rest of the claim
we need to estimate $|y_i-\tilde{y_i}|$ and $|\mu_i-\tilde{\mu_i}|$. 
To this end we observe that by construction
\begin{align*}
o_{A_0}(1) = \|z_{\mu_i,y_i}- z_{\tilde{\mu_i},\tilde{y_i}}\|_{\Did} 
= \|z_{\tilde{\mu_i}/\mu_i,\tilde{y_i}-y_i}-z_{1,0}\|_{\Did}. 
\end{align*}
Since
\begin{align*}
\lim_{\mu+\mu^{-1}+|y|\to \infty}\|z_{\mu,y}-z_{1,0}\|_{\Did}^2= 2\|z_{1,0}\|_{\Did}^2  
\end{align*}
we see that 
there is $R_1=R_1(A_0)>0$ such that
\begin{align*}
(R_1)^{-1} \le \tilde{\mu_i}/\mu_i \le R_1 \text{ and }  |y_i-\tilde{y_i}| \le R_1.
\end{align*}
Now by explicit calculations or elliptic regularity \cite{BrezisKato79} we have
\begin{align*}
\max\Big\{ |y_i-\tilde{y_i}|,\Big|\frac{\tilde{\mu_i}}{\mu_i}-1\Big|\Big\} \le {\rm const}(A_0) 
\|z_{\frac{\tilde{\mu_i}}{\mu_i},\tilde{y_i}-y_i}-z_{1,0}\|_{\Did} = o_{A_0}(1),
\end{align*}
which gives the claim.
\end{proof}
\section{Expansion of the perturbation terms $w$ and $\vec{\a}$}
\label{s:expansion}
For the rest of the paper we will only treat the case $N=3$. Unless otherwise indicated, integration
extends over $\rz^3$ and is done with respect to the variable $x$. Moreover, we will write $k$ instead of
$k_\theta$ when there is no possibility of confusion to avoid cumbrous subindexing .\\
From the change of coordinates $x \mapsto \mu x+y$, H\"older's and Sobolev's inequality we get
\begin{lemma}
\label{l:estimates1}
Let $y \in \rz^3$, $\tau>0$ and 
$f,r: \rz^3 \to \rz$ measurable such that 
\begin{align*}
|r(x)| \le C_r |x-y|^{\sigma} \ \text{ in } B_1(y), \
|r(x)| \le C_r |x-y|^{\tilde \sigma}\ \text{ in } \rz^3 \setminus B_1(y),\\
|f(x)| \le C_r |x|^{-s} \text{ in }B_1(0), \ |f(x)| \le C_r |x|^{-m} \text{ in }\rz^3\setminus B_1(0), 
\end{align*}
for some $C_r,m,s >0$ and $0\le \tilde{\sigma},\sigma$. Then there is $C=C(\tau,C_r)>0$ such that 
for $v,v_1,v_2 \in \Did$:\\
For $ 0\le \tilde{\sigma},\sigma\le m-3-\tau$ and $s+\tau \le 3+\sigma$ there holds 
\begin{align*}
\big| \int r(x) \mu^{-3} f\Big(\frac{x-y}{\mu}\Big)\big| 
\le C 
\big(\charak_{(0,1]}(\mu)\mu^\sigma+\charak_{(1,\infty)}(\mu)(\mu^{s-3}+\mu^{\tilde{\sigma}})\big),
\end{align*}
if $ 0\le \tilde{\sigma},\sigma\le m-\frac{5}{2}-\tau$ and $s+\tau <\sigma +\frac{5}{2}$ then 
\begin{align*}
\big\| \int \frac{r(x)}{\mu^{\frac52}} f\Big(\frac{x-y}{\mu}\Big) \cdot \big\|_{\Did} 
\le C 
\big(\charak_{(0,1]}(\mu)\mu^\sigma+\charak_{(1,\infty)}(\mu)(\mu^{s-\frac{5}{2}}+\mu^{\tilde{\sigma}})\big),
\end{align*}
if $0 \le \tilde{\sigma}, \sigma \le m-2-\tau$ and $s+\tau <\sigma +2$ then  
\begin{align*}
\sup_{\|v\|\le 1}
\big\| \int r(x) \mu^{-2} &f\Big(\frac{x-y}{\mu}\Big) v \cdot \big\|_{\Did}\\
&\le C 
\big(\charak_{(0,1]}(\mu)\mu^\sigma+\charak_{(1,\infty)}(\mu)(\mu^{s-2}+\mu^{\tilde{\sigma}})\big),
\end{align*}
if $0\le \tilde{\sigma}, \sigma \le m-\frac{3}{2}-\tau$ and $s+\tau <\sigma +\frac{3}{2}$ then we have  
\begin{align*}
\begin{split}
\sup_{\|v_1\|,\|v_2\|\le 1}
\big\| \int r(x) \mu^{-\frac32} &f\Big(\frac{x-y}{\mu}\Big) v_1 v_2 \cdot \big\|_{\Did} \\
&\le C 
\big(\charak_{(0,1]}(\mu)\mu^\sigma+\charak_{(1,\infty)}(\mu)(\mu^{s-\frac{3}{2}}+\mu^{\tilde{\sigma}})\big).    
\end{split}
\end{align*}
\end{lemma}
Using the above estimates we may prove the main ingredient for the 
finite dimensional reduction.
\begin{lemma}
\label{p:implicit} 
Suppose $k \in C^5(\rz^3)$ and there are $A_0,B_0,B_1>0$ such that
\begin{align*}
\max\big(B_0,B_1,\sup_{|m|\le 5}\|D^{m}k\|_\infty\big) \le A_0 \text{ and }\\
A_0^{-1} \le 1+t k(x) \; \forall (x,t) \in \rz^3 \times [-B_0,B_1].
\end{align*}
Then there exist  $\rho_0=\rho_0(A_0)>0$, $t_0= t_0(A_0)>0$, 
an upper continuous function $\mu_0:\rz \to \rz_{+} \cup \{\infty\}$, 
depending only on $A_0$, 
and two functions $w:\Omega \to \Did$ and $\vec\alpha:\Omega \to \rz^{4}$, where
\begin{align*}
\Omega:= \{(t,\mu,y) \in [-B_0,B_1] \times (0,+\infty)\times \rz^3 \where 0<\mu<\mu_0(t)\},\\
\mu_0(t) = +\infty \ \text{ if } |t|\le t_0,  
\end{align*}
such that for any $(t,\mu,y) \in \Omega$
\begin{align}
\label{eq:12}
w(t,\mu,y)\ \text{ is orthogonal to }\ T_{z_{\mu,y}}Z\\
\label{eq:13}
f_{t}'\big(z_{\mu,y}+w(t,\mu,y)\big)= \vec \alpha(t,\mu,y) \cdot \dot{\xi}_{\mu,y} \in T_{z_{\mu,y}}Z\\
\label{eq:9}
\|w(t,\mu,y) -  w_0(t,\mu,y)\|+ \|\vec{\alpha}(t,\mu,y)\| < \rho_0,
\end{align}
where $\{{\dot\xi^i}_{\mu,y} \where i=0 \dots 3\}$ denotes the orthonormal basis of $T_{z_{\mu,y}}Z$ 
given in \eqref{eq:11} and
\begin{align*}
w_0(t,\mu,y) := \big((1+t k(y))^{-\frac{1}{4}}-1\big)z_{\mu,y}.
\end{align*}
The functions $w$ and $\vec \alpha$ are of class $C^{2}$ 
and unique in the sense that if $(v,\vec \beta)$ satisfies
\eqref{eq:12}-\eqref{eq:9} for some $(t,\mu,y)\in \Omega$ then $(v,\vec \beta)$ is given by
$(w(t,\mu,y),\vec \alpha(t,\mu,y))$. Moreover, we have for $1 \le j \le 3$
\begin{align}
\label{eq:14}
\begin{split}
&\|w(t,\mu,y) - \sum_{i=0}^2 w_i(t,\mu,y)\| 
+\|\vec{\alpha}(t,\mu,y)_0-\sum_{i=1}^2 \vec{\alpha}_i(t,\mu,y)_0\|\\ 
&\hspace{7em} \leq O_{A_0}\Big(t\big(|\nabla k(y)|^2\min(1,\mu^{2})+\min(1,\mu^{\frac94})\big)\Big),\\
&\|\vec{\alpha}(t,\mu,y)_j-\sum_{i=1}^2 \vec{\alpha}_i(t,\mu,y)_j\| \leq  
O_{A_0}\big(t \min(1,\mu^{\frac94})\big),  
\end{split}  
\end{align}
where 
\begin{align*}
\vec{\alpha}_1(t,\mu,y):= -t \min(1,\mu) (1+t k(y))^{-\frac{5}{4}}\,
\frac{\pi}{3^{\frac14}\sqrt{5}}
\binom{0}{\nabla k(y)},\\
\vec{\alpha}_{2}(t,\mu,y) := -t \min(1,\mu^2) (1+t k(y))^{-\frac{5}{4}}
\frac{\pi}{3^{\frac14}\sqrt{5}} \binom{\laplace k(y)}{\vec 0}.
\end{align*}
and 
\begin{align*}
w_1(t,\mu,y) &:=  t \min(1,\mu) (1+t k(y))^{-\frac{5}{4}} \ 
{\mathcal T}_{y}\circ {\mathcal U}_\mu \big({\tilde w}_1(y)\big),\\
{\tilde w}_1(y) &:= {\mathcal F_0^{-1}}\big(\int \nabla k(y) x (z_{1,0})^{5}\cdot \big),\\
w_2(t,\mu,y) &:= t \min(1,\mu^2) (1+t k(y))^{-\frac{5}{4}} \ 
{\mathcal T}_{y}\circ {\mathcal U}_\mu \big({\tilde w}_2(y)\big),\\
{\tilde w}_2(y) &:=  {\mathcal F_0^{-1}} 
\Big(
\frac{1}{2} \int D^2k(y)x^2 (z_{1,0})^{5} \cdot 
\Big),
\end{align*}
The operator ${\mathcal F_0^{-1}} \in {\mathcal{L}}(\Did, T_{z_{1,0}}Z^\perp)$ is defined by
\begin{align*}
{\mathcal F_0^{-1}} := \big(f_0''(z_{1,0})|_{T_{z_{1,0}}Z^\perp}\big)^{-1} \circ {\text{Proj}}_{T_{z_{1,0}}Z^\perp}. 
\end{align*}
\end{lemma}
\begin{proof}
Define $H: \rz \times (0,\infty)\times \rz^3 \times\Did \times \rz^{4} \to \Did \times \rz^{4}$
\begin{align*}
H(t,\mu,y,w,\vec \alpha):= \big(f_t'(z_{\mu,y}+w)-\vec \alpha \cdot \dot{\xi}_{\mu,y},
(\langle w,(\dot{\xi}_{\mu,y})_l\rangle)_l\big).  
\end{align*}
If $H(t,\mu,y,w,\vec\alpha)=(0,0)$ then $w$ satisfies \eqref{eq:12}-\eqref{eq:13}. 
We have
\begin{align}
\label{eq:91}
\bigg(\frac{\rand H}{\rand(w,\vec\alpha)}(t,\mu,y,w,\vec\a)\bigg)\binom{\phi}{\vec\beta} =
\big(f_t''(z_{\mu,y}+w)\phi - \vec\beta \dot{\xi}_{\mu,y},(\langle\phi,(\dot{\xi}_{\mu,y})_l\rangle)_l\big).  
\end{align}
Note that 
\begin{align}
\label{eq:5}
\begin{split}
\bigg\langle \bigg(\frac{\rand H}{\rand(w,\vec\alpha)}(0,\mu,y,0,0)\bigg)&\binom{w}{\vec\beta},
\big(f_0''(z_{\mu,y})w-\vec\beta \cdot \dot{\xi}_{\mu,y},(\langle w,(\dot{\xi}_{\mu,y})_l\rangle)_l \big)\bigg\rangle\\
&= \|f_0''(z_{\mu})w\|^2 + |\vec\beta|^2+ |(w,\dot{\xi}_{\mu,y})_i|^2. 
\end{split}
\end{align}
From \eqref{eq:3} and \eqref{eq:5} we infer that $\big(\frac{\rand H}{\rand(w,\vec\alpha)}
(0,\mu,y,0,0)\big)$ is an injective Fredholm operator of index zero, hence invertible and
\begin{align}
\label{eq:4}
\bigg\|\bigg(\frac{\rand H}{\rand(w,\vec\alpha)}(0,\mu,y,0,0)\bigg)^{-1}\bigg\| \le
1+\|(f_0''(z_{\mu,y}))^{-1}\|_{{\mathcal{L}}(T_{z_{\mu,y}}Z^\perp)} =: C_*.   
\end{align}
Clearly, $H(t,\mu,y,w,\vec \alpha)=(0,0)$ if and only if $(w,\vec\alpha)= F_{t,\mu,y}(w,\vec\alpha)$, where 
\begin{align*}
F_{t,\mu,y}(w,\vec\alpha):= -\bigg(\frac{\rand H}{\rand(w,\vec\alpha)}(0,\mu,y,0,0)\bigg)^{-1}
H(t,\mu,y,w,\vec\alpha) + (w,\vec\alpha).
\end{align*}
We will prove that $F_{t,\mu,y}(w,\vec\alpha)$ is a contraction in
some ball 
\[B_\rho\bigg(\sum_{i=0}^2w_i(t,\mu,y),\sum_{i=1}^2 {\vec \alpha}_i(t,\mu,y)\bigg)\]
for any radius $\rho$ such that
\[O_{A_0}\Big(t\big(|\nabla k(y)|\min(1,\mu^{2})+\min(1,\mu^{2+\frac14})\big)\Big) \le \rho \le \rho_0,\]
where $\rho_0=\rho_0(A_0)$ will be chosen later.\\
To this end we fix $\rho>0$ and $(w,\vec\alpha)\in B_\rho(0,0)$. 
In the sequel we will suppress the dependence of $\vec{\alpha}_i$ and $w_i$ on 
$t$, $\mu$ and $y$.
From \eqref{eq:4} and Sobolev's inequality 
\begin{align}
\frac{1}{C_*}\|&F_{t,\mu,y}\big(w+\sum_{i=0}^2w_i,\vec\alpha+\sum_{i=1}^2 {\vec \alpha}_i\big)
-\big(\sum_{i=0}^2w_i,\sum_{i=1}^2{\vec \alpha}_i\big)\| \notag \\
&\le \|f_t'\big(z_{\mu,y}+w+\sum_{i=0}^2w_i\big)
- \sum_{i=1}^2{\vec \alpha}_i \cdot \dot{\xi}_{\mu,y} -f_0''(z_{\mu,y})w\|\notag \\
&\le \|\langle z_{\mu,y}+\sum_{i=0}^2w_i - \sum_{i=1}^2{\vec \alpha}_i \cdot \dot{\xi}_{\mu,y}, \cdot \rangle
 -5 \int (z_{\mu,y})^{4}w \cdot \notag \\
\label{eq:6}
&\quad -\int(1+ t k(x))\Big( \sum_{l=0}^2 \binom{5}{l} (z_{\mu,y}+w_0)^{5-l}(w_1+w_2+w)^l\Big) \cdot \| \notag \\
&\quad +O_{A_0}(\|w_1+w_2+w\|^{3}).
\end{align}
Obviously, $(1+ t k(y)) (z_{\mu,y}+w_0)^{4} = (z_{\mu,y})^{4}$ and
\begin{align}
\label{eq:17}
\langle z_{\mu,y}+w_0,\varphi\rangle  =  \int (1+ t k(y)) (z_{\mu,y}+w_0)^{5}\varphi.
\end{align}
Inserting this in \eqref{eq:6} and using Lemma \ref{l:estimates1} we get
\begin{align}
\label{eq:10}
\frac{1}{C_*}\|F_{t,\mu,y}&\big(w+\sum_{i=0}^2w_i,\vec\alpha+\sum_{i=1}^2 {\vec \alpha}_i\big)
-\big(\sum_{i=0}^2w_i,\sum_{i=1}^2{\vec \alpha}_i\big)\| \notag \\
&\le \|f_0''(z_{\mu,y})(w_1+w_2)
-\langle({\vec \alpha}_1+ {\vec \alpha}_2)\cdot \dot{\xi}_{\mu,y}, \cdot \rangle \notag \\
&\quad - t \int (k(x)-k(y)) (z_{\mu,y}+w_0)^{5} \cdot \|\notag \\
&\quad + O_{A_0}\Big(t^2\big(|\nabla k(y)|\min(1,\mu^2)+\min(1,\mu^3)\big)\Big)\notag \\
&\quad + O_{A_0}\Big(t\min(1,\mu)\|w\|+\|w\|^2\Big).  
\end{align}
From \eqref{eq:7}-\eqref{eq:2} and the definition of $w_1$ and $w_2$ we infer 
\begin{align*}
f_0''(z_{\mu,y})w_1 &= 
t \min(1,\mu^{-1}) {\text{Proj}}_{T_{z_{\mu,y}}Z^\perp} 
\big(\int \nabla k(y)(x-y) (z_{\mu,y}+w_0)^{5}\cdot \big),\\
f_0''(z_{\mu,y})w_2 &= 
t \min(1,\mu^{-2}) {\text{Proj}}_{T_{z_{\mu,y}}Z^\perp} 
\Big(\int \frac12 D^2k(y) (x-y)^2 (z_{\mu,y}+w_0)^{5} \cdot \Big).
\end{align*}
To find $\vec{\a}_1$ and $\vec{\a}_2$ we observe that since $(\dot{\xi}_{\mu,y})_0$ is even
and $(\dot{\xi}_{\mu,y})_i$ is odd for $1 \le i \le 3$ we have
\begin{align*}
\int \nabla k(y)(x-y) (z_{\mu,y}+w_0)^{5} (\dot{\xi}_{\mu,y})_0 = 0,\\
\int \frac12 D^2k(y) (x-y)^2 (z_{\mu,y}+w_0)^{5} (\dot{\xi}_{\mu,y})_i &= 0 \text{ for } 1\le i \le 3.
\end{align*}
For $1 \le i \le 3$ we get from \eqref{eq:A44}-~\eqref{eq:A8}
\begin{align*}
\int \nabla k(y)(x-y) (z_{\mu,y}+w_0)^{5} (\dot{\xi}_{\mu,y})_i 
&= \mu\, \frac{\pi}{3^{\frac{1}{4}}\sqrt{5}} (1+t k(y))^{-\frac{5}{4}} \frac{\rand k}{\rand x_i}(y). 
\end{align*}
For $i=0$ we get
\begin{align*}
\int \frac12 D^2k(y) (x-y)^2 (z_{\mu,y}+w_0)^{5} (\dot{\xi}_{\mu,y})_i
= \mu^2 \frac{\pi}{3^{\frac{1}{4}}\sqrt{5}} (1+t k(y))^{-\frac54}  \laplace k(y).
\end{align*}
Finally, we obtain
\begin{align}
\label{eq:66}
\begin{split}
f_0''&(z_{\mu,y})(w_1+w_2) - (\vec{\alpha}_1+\vec{\alpha}_2)\cdot \dot{\xi}_{\mu,y}\\
&= t \Big(\int \sum_{\ell =1}^2 \frac{D^\ell k(y)}{\ell!}(x-y)^\ell (z_{\mu,y}+w_0)^{5}\cdot\Big)
+O_{A_0}\big(t \charak_{(1,\infty)}(\mu)\big),    
\end{split}
\end{align}
which implies together with \eqref{eq:10} and Lemma \ref{l:estimates1}
\begin{align}
\label{eq:18}
\frac{1}{C_*}\|&F_{t,\mu,y}\big(w+\sum_{i=0}^2w_i,\vec\alpha+\sum_{i=1}^2 {\vec \alpha}_i\big)
-\big(\sum_{i=0}^2w_i,\sum_{i=1}^2{\vec \alpha}_i\big)\| \notag \\
&\le O_{A_0}\Big(t\big(|\nabla k(y)|\min(1,\mu^2)+\min(1,\mu^{2+\frac14})\big)\Big) \notag\\
&\quad +O_{A_0}\Big(t\min(1,\mu)\|w\|+\|w\|^2\Big).
\end{align}
Consequently, if we fix $0<\rho_0<1/4$ we obtain functions $\mu_0$ and $t_0$ depending on $\rho_0$ and $A_0$ such that
$F_{t,\mu,y}$ maps $B_\rho(\sum_{i=0}^2 w_i,\sum_{i=1}^2 \vec{\a}_i)$
into itself for every $(t,\mu,y) \in \Omega$ and $\rho>0$ satisfying
\begin{align}
\label{eq:22}
{\rm const}(A_0,\rho_0)|t|\big(|\nabla k(y)|\min(1,\mu^2)+\min(1,\mu^{2+\frac14})\big) <\rho\le \rho_0. 
\end{align}
To show that $F_{t,\mu,y}$ is a contraction we fix $\rho>0$ and two vectors  
$(v_1,\vec{\beta}_1)$ and $(v_2,\vec{\beta}_2)$ in $B_\rho(0,0)$.
Then using Lemma \ref{l:estimates1} and \eqref{eq:4}
\begin{align*}
&\frac{\|F_{t,\mu,y}(\sumg_{i=0}^2 w_i+v_1,\sumg_{i=1}^2 \vec{\a}_i+\vec{\beta}_1)
-F_{t,\mu,y}(\sumg_{i=0}^2 w_i+v_2,\sumg_{i=1}^2 \vec{\a}_i+\vec{\beta}_2)|}
{C_*\|(v_1,\beta_1)-(v_2,\beta_2)\|} \\
&\quad \le \int_0^1 \|f_t''\big(z_{\mu,y}+\sum_{i=0}^2 w_i+v_1+s(v_2-v_1)\big)-f_0''(z_{\mu,y})\| \di s\\
&\quad \le \int_0^1 \|f_t''\big(z_{\mu,y}+w_0\big)-f_0''(z_{\mu,y})\| \di s +O_{A_0}(\rho+\|w_1+w_2\|)\\
&\quad \le O_{A_0}\big(\rho+t\min(1,\mu)\big). 
\end{align*}
Thus we have shown that there are $\rho_0>0$ and two function $\mu_0$ and $t_0$ depending only on $A_0$, as claimed above, such that 
$F_{t,\mu,y}$ is a contraction in $B_\rho(\sum_{i=0}^2 w_i,\sum_{i=1}^2 \vec{\a}_i)$ 
for every $(t,\mu,y)\in \Omega$ and $\rho>0$ satisfying \eqref{eq:22}.\\  
From Banach's fixed-point theorem we deduce the existence and uniqueness of the functions $w$ and $\vec{\a}$.
The usual inverse function theorem yields the $C^2$ dependence. The estimates in \eqref{eq:14} hold due to the uniqueness
of the fixed-point and because $F_{t,\mu,y}$ is a contraction for every $\rho>0$ satisfying \eqref{eq:22}.
\end{proof}
We need a precise expansion of $w$ and $\vec\alpha$ in critical points $y$ of $k$ 
in terms of $\mu$ and $t$.
This will be done up to order $5$ in $\mu$. We later see that we may assume $|\nabla k(y)|$ to be of order
$O(\mu^2)$. First we compute ${\tilde w}_{2}(y)$ in terms of the eigenfunctions of $f_0''(z_{1,0})$.
\begin{lemma}
\label{l:comp:w2}
Under the assumptions of Lemma \ref{p:implicit} we have  
\begin{align*}
{\tilde w}_{2}(y) &= \frac{3^\frac54 \sqrt{\pi}}{8}
\sum_{j=0}^{\infty} \frac{\Gamma(j+\frac{9}{2}) (5+2j) a_{2,j}}{\Gamma(j+6)}\\ 
&\quad \qquad \times \bigg(\sum_{1\le l<m\le 3} \psi_{l,m}(y) \Phi^{1,0}_{2,j,(l,m)}
+ \sum_{l=2}^3 \psi_{l,l}(y) \Phi^{1,0}_{2,j,(l,l)}\bigg) \\
&\quad -\frac{\pi 3^\frac34}{16} \laplace k(y) \Phi^{1,0}_{0,0,1}
+ \frac{3^\frac14 \pi}{4}\laplace k(y)
\sum_{j=2}^{\infty} \frac{\Gamma(j+\frac52)(1+2j)a_{0,j}}{\Gamma(j+2)(j+3)(j-1)}  \Phi^{1,0}_{0,j,1},     
\end{align*}
where we use the basis defined in \eqref{eq:37}, \eqref{eq:115} and 
\begin{align*}
\psi_{l,m}(y) &:= \frac{2 \sqrt{\pi}}{\sqrt{15}}  \frac{\rand^2k(y)}{\rand x_l \rand x_m},
&\psi_{2,2}(y) &:= \frac{\sqrt{\pi}}{\sqrt{15}}
\Big( \frac{\rand^2k(y)}{\rand x_2^2} -\frac{\rand^2k(y)}{\rand x_1^2}\Big),\\
\psi_{3,3}(y) &:= \frac{\sqrt{\pi}}{\sqrt{5}}
\Big( \frac{\rand^2k(y)}{\rand x_3^2} -\frac{\laplace k(y)}{3}\Big).
\end{align*}
\end{lemma}
\begin{proof}
To prove the claim we observe that
if 
\begin{align*}
{\text{Proj}}_{T_{z_{1,0}}Z^\perp}\big(\frac{1}{2} \int D^2k(y)(x)^2 (z_{1,0})^{2^*-1} \cdot \big)
= \sum_{i+j\neq 1} \sum_{l=1}^{c_i} \beta_{i,j,l} \Phi^{1,0}_{i,j,l}, 
\end{align*}
then \eqref{eq:115} implies
\begin{align}
\label{eq:42}
{\tilde w}_{2}(y) = \sum_{i+j\neq 1} \sum_{l=1}^{c_i} \frac{\beta_{i,j,l}}{\lambda_{i,j}} \Phi^{1,0}_{i,j,l}.  
\end{align}
To this end we note that in the basis given in \eqref{eq:37}
\begin{align}
\label{eq:59}
\begin{split}
\frac{1}{2} D^2k(y)(x)^2
&= \frac{\laplace k(y)}{6} \|x\|^2 
+ \sum_{l=2}^3 \psi_{l,l}(y) Y_{2,(l,l)}(x)\\
&\qquad +\sum_{1 \le l<m \le 3} \psi_{l,m}(y) Y_{2,(l,m)}(x).  
\end{split}
\end{align}
Consequently,
\begin{align*}
{\tilde w}_{2}(y) \in \langle \Phi^{1,0}_{2,j,l},\Phi^{1,0}_{0,j,1} \where j \in \nz_0,\, 1\le l \le c_2 \rangle.   
\end{align*}
Using \eqref{eq:A1}-\eqref{eq:A8} we find for $n=(l,l)$ or $n=(l,m)$
\begin{align*}
\int \psi_n Y_{2,n}(x/|x|) |x|^2 (z_{1,0})^{5} \Phi^{1,0}_{2,j,n}
= \psi_n a_{2,j} 3^{\frac{5}{4}} \frac{(j+1) \Gamma(3/2) \Gamma(j+7/2)}{2 \Gamma(j+5)}\\ 
\int |x|^2 z_{1,0}^{5} \Phi^{1,0}_{0,j,0}
= a_{0,j} 3^{\frac{5}{4}} 2 \sqrt{\pi} 
\bigg( \frac{\Gamma(\frac{1}{2}) \Gamma(j+\frac{3}{2})}{2\Gamma(j+2)}
-\delta_{0,j} \frac{\Gamma(\frac{3}{2})^2}{2\Gamma(3)}\bigg).
\end{align*}
Now, the claim follows from \eqref{eq:42} and \eqref{eq:59}.
\end{proof}

\begin{remark}
\label{r:global}
Under the assumptions of Lemma \ref{p:implicit} 
an explicit calculation together with \eqref{eq:14} yields for $(t,\mu,y) \in \Omega$, where we suppress the
dependence of $w$ and $\vec{\a}$ on $(t,\mu,y)$,
\begin{align}
f_t''(z_{\mu,y}+w) \phi &= f_0''(z_{\mu,y}) \phi -5t \int (k(x)-k(y)) (z_{\mu,y}+w_0)^4 \phi \cdot \notag\\
&\quad -20 \int (1+t k(x))(z_{\mu,y}+w_0)^3(w-w_0) \phi \cdot\notag\\ 
&\quad -30 \int (1+t k(x))(z_{\mu,y}+w_0)^2(w-w_0)^2 \phi \cdot\notag\\
\label{eq:77}
&\quad + t^3 O_{A_0}(|\nabla k(y)|^3 \min(1,\mu^3)+\min(1,\mu^6)),
\end{align}
which implies by Lemma \ref{l:estimates1}
\begin{align*}
\|f_t''(z_{\mu,y}+w)- f_0''(z_{\mu,y})\| \le |t| O_{A_0}(|\nabla k(y)| \min(1,\mu)+\min(1,\mu^\frac74)). 
\end{align*}
Consequently, from \eqref{eq:91},
after decreasing $\mu_0$ and $t_0$ if necessary, 
\begin{align*}
\Big\| \frac{\rand H}{\rand(w,\vec\alpha)}(t,\mu,y,w,\vec{\a})- 
\frac{\rand H}{\rand(w,\vec\alpha)}(0,\mu,y,0,0) \Big\| \le \frac12 C_*
\end{align*}
Thus we may assume $\frac{\rand H}{\rand(w,\vec\alpha)}(t,\mu,y,w,\vec{\a})$ is invertible and its inverse
is uniformly bounded with respect to $(t,\mu,y)\in \Omega$. 
\end{remark}
We begin the expansion of $\vec{\a}$ by computing the third order term.
\begin{lemma}
\label{l:expansion:3}
Under the assumptions of Lemma \ref{p:implicit} we have as $\mu \to 0$
\begin{align*}
\|\vec{\alpha}(t,\mu,y)-\sum_{j=1}^3\vec{\alpha}_{j}(t,\mu,y)\| 
= t O_{A_0}(\mu^{3+\frac12} +|\nabla k(y)|^2\mu^2),            
\end{align*}
where $\vec{\alpha}_{1}$, $\vec{\alpha}_{2}$ are defined in Lemma \ref{p:implicit} 
and  $\vec{\alpha}_{3}$ is given by
\begin{align*}
\vec{\alpha}_{3}(t,\mu,y)_i &:= -t \mu^3 \frac{3^{\frac14} \pi}{2 \sqrt{15}} 
(1+t k(y))^{-\frac{5}{4}} \frac{\rand}{\rand x_i}\laplace k(y),
\end{align*}
for $i=1 \dots 3$ and
\begin{align*}
\vec{\alpha}_{3}(t,\mu,y)_0 &:= -t \mu^3 (1+t k(y))^{-\frac{5}{4}} \frac{3^\frac34 4}{\pi \sqrt{5}}
\cpvint{}{} \big(k(x+y)- T_{k(\cdot+y),0}^3(x)\big)\frac{1}{|x|^{6}}.
\end{align*}
\end{lemma}
\begin{proof}
In the sequel we will suppress the dependence of $w$ and $\vec{\alpha}$ on $t$, $\mu$ and $y$, when there is no 
possibility of confusion. Moreover, we always assume $0<\mu \le 1$.\\ 
As in Lemma \ref{p:implicit} we infer from Lemma \ref{l:estimates1}, \eqref{eq:17}, and the
definition of $w_1$
\begin{align}
\label{eq:72}
f_t'(z_{\mu,y}+w)  
&= f_0''(z_{\mu,y})(w-w_0)
-t \int \big(k(x)-k(y)\big)\,(z_{\mu,y}+w_0)^{5} \cdot \notag \\ 
&\quad-5 t \int \big(k(x)-T_{k,y}^1(x)\big)\,(z_{\mu,y}+w_0)^{4} (w-w_0-w_1) \cdot\notag \\
&\quad- 10 \int (1+t k(x)) (z_{\mu,y}+w_0)^{3} (w-w_0-w_1)^2 \cdot \notag \\
&\quad + t^2 O_{A_0}(|\nabla k(y)|^2\mu^2+|\nabla k(y)|\mu^3+\mu^6).
\end{align}
As $f_t'(z_{\mu,y}+w)-\vec{\alpha} \dot\xi_{\mu,y}=0$ and $(\dot\xi_{\mu,y})_i \in N\big( f_0''(z_{\mu,y})\big)$
we obtain from \eqref{eq:66} and testing \eqref{eq:72} with $(\dot\xi_{\mu,y})_i$ 
\begin{align}
\label{eq:114}
\big(\vec{\alpha}-\sum_{j=1}^2\vec{\alpha}_j\big)_i &= -t
\int \big(k(x)-T_{k,y}^2(x)\big)\,(z_{\mu,y}+w_0)^{5}(\dot \xi_{\mu,y})_i \notag \\ 
&\quad-5 t \int \big(k(x)-T_{k,y}^1(x)\big)\,(z_{\mu,y}+w_0)^{4} (w-w_0-w_1) (\dot \xi_{\mu,y})_i\notag \\
&\quad- 10 \int (1+t k(x)) (z_{\mu,y}+w_0)^{3} (w-w_0-w_1)^2 (\dot \xi_{\mu,y})_i \notag \\
&\quad + t^2 O_{A_0}(|\nabla k(y)|^2\mu^2+|\nabla k(y)|\mu^3+\mu^6).
\end{align}
By Lemma \ref{l:estimates1} and \ref{p:implicit} we have for $0\le i\le 3$
\begin{align}
\label{eq:129}
\begin{split}
\int \big(k(x)-T_{k,y}^1(x)\big)\,(z_{\mu,y}+w_0)^{4} (w-w_0-w_1) (\dot \xi_{\mu,y})_i =  O_{A_0}(t\mu^{3+\frac12}),\\
\int (1+t k(x)) (z_{\mu,y}+w_0)^{3} (w-w_0-w_1)^2 (\dot \xi_{\mu,y})_i =  O_{A_0}(t\mu^{4}).     
\end{split}
\end{align}
If $1 \le i \le 3$ then we obtain from Lemma \ref{l:estimates1}
\begin{align*}
\int \big(k(x)-&T_{k,y}^2(x)\big) \,(z_{\mu,y}+w_0)^{5}(\dot \xi_{\mu,y})_i\\
&= \int \frac{1}{6}D^3 k(y)(x-y)^3\,(z_{\mu,y}+w_0)^{5}(\dot \xi_{\mu,y})_i + O_{A_0}(\mu^{3+\frac12}),
\end{align*}
and from \eqref{eq:A70}
\begin{align}
\label{eq:100}
\int \frac{1}{6}&D^3 k(y)(x-y)^3\,(z_{\mu,y}+w_0)^{5}(\dot \xi_{\mu,y})_i
= \mu^3\frac{3^\frac14 \pi}{2 \sqrt{15}}(1+t k(y))^{-\frac54} \frac{\rand}{\rand x_i}\laplace k(y). 
\end{align} 
Hence, the assertion of the lemma for $1 \le i \le 3$ follows from \eqref{eq:114}-\eqref{eq:100}.\\
To treat the remaining case $i=0$ we use the fact that $D^3 k(y)(x)^3$ is odd and get
\begin{align}
\label{eq:32}
\begin{split}
\intg_{B_1(y)} \big(k(x)-&T^2_{k,y}(x)\big)\,(z_{\mu,y}+w_0)^{5}(\dot \xi_{\mu,y})_0\\
&=\intg_{B_1(y)} \big(k(x)-T^3_{k,y}(x)\big)\,(z_{\mu,y}+w_0)^{5}(\dot \xi_{\mu,y})_0  
\end{split}
\end{align}
Since, by Lemma \ref{l:estimates1}, 
\begin{align*}
\int_{\rz^3}\big|k(x)-T^3_{k,y}(x)\big| \,\mu^{-3}(1+\mu^{-2}|x-y|^2)^{-4}
= O_{A_0}(\mu^4), 
\end{align*}
there holds after a translation $x \to x+y$
\begin{align*}
\intg_{B_1(y)} &\big(k(x)-T^3_{k,y}(x)\big)\,(z_{\mu,y}+w_0)^{5}(\dot \xi_{\mu,y})_0\\
&= \frac{3^{\frac{3}{4}} 4 \mu^3}{\pi \sqrt{5}(1+t k(y))^{\frac{5}{4}}}
\intg_{B_1(0)} \frac{\big(k(x+y)-T^3_{k(\cdot +y),0}(x)\big)}{(\mu^2+|x|^2)^{3}}+ O_{A_0}(\mu^4).
\end{align*}
From
\begin{align*}
(\mu^2+|x|^2)^{-3}-|x|^{-6} &= -(\mu^2+|x|^2)^{-3}\big(3|x|^{-2}\mu^2+3|x|^{-4}\mu^4+|x|^{-6}\mu^6 \big)
\end{align*}
and Lemma \ref{l:estimates1}, we infer
\begin{align}
\label{eq:74}
\int_{\rz^3} &\big|k(x+y)-T^3_{k(\cdot +y),0}(x)\big|
\big|(\mu^2+|x|^2)^{3}-|x|^{-6}\big| = O_{A_0}(\mu).
\end{align}
Hence,
\begin{align*}
\int \big(&k(x)-T^2_{k,y}(x)\big)\,(z_{\mu,y}+w_0)^{5}(\dot \xi_{\mu,y})_0 \\
&= \frac{3^{\frac{3}{4}}4\mu^3}{\pi\sqrt{5} (1+t k(y))^{\frac{5}{4}}} 
\Bigg( \int_{B_1(0)} \big(k(x+y)-T^3_{k(\cdot +y),0}(x)\big)\,|x|^{-6}+\\
&\qquad \int_{\rz^3\setminus B_1(0)} \big(k(x+y)-T^2_{k(\cdot +y),0}(x)\big)\,|x|^{-6}\Bigg)
+ O_{A_0}(\mu^4)\\
&= \frac{3^\frac34 4}{\pi \sqrt{5}} \mu^3 (1+t k(y))^{-\frac{5}{4}} 
\cpvint{}{} \Big(k(x+y)- T^3_{k(\cdot +y),0}(x)\Big)\frac{1}{|x|^{6}}+ O_{A_0}(\mu^4),
\end{align*}
which ends the proof.
\end{proof}

\begin{remark}
\label{r:1_t_alpha_continuous}
From Lemma \ref{p:implicit}, \eqref{eq:72}, and \eqref{eq:114} we see that
\begin{align*}
\frac{1}{t \mu} \vec{\alpha}(t,\mu,y) 
\end{align*}
is a well defined, continuous function for $(t,\mu,y) \in \Omega$.
\end{remark}

\begin{lemma}
\label{l:expansion:4}
Under the assumptions of Lemma \ref{p:implicit} we have as $\mu \to 0$
\begin{align*}
|\vec{\alpha}(&t,\mu,y)-\sum_{j=1}^4\vec{\alpha}_{j}(t,\mu,y)| 
= tO_{A_0}(\mu^{4+\frac12})\\
&+t^2 O_{A_0}\big(\mu^2|\nabla k(y)|^2+\mu^3|\nabla k(y)|+ \mu^4 |\laplace k(y)|+ \mu^{4+\frac14}\big),            
\end{align*}
where for $1\le i \le 3$
\begin{align*}
\vec{\alpha}_{4}(t,\mu,y)_i &=  - t \mu^4 (1+t k(y))^{-\frac{5}{4}} \frac{3^\frac34 8}{\pi \sqrt{5}} 
\cpvint{}{} \Big(k(x+y)-  T^3_{k(\cdot +y),0}(x)\Big)\frac{x_i}{|x|^{8}}  
\end{align*}
and
\begin{align*}
\vec{\alpha}_{4}(t,\mu,y)_0 &= t \mu^4  (1+t k(y))^{-\frac{5}{4}} 
\frac{3^\frac34 \pi \sqrt{5}}{30} \laplace^2 k(y)\\ 
&\quad\,- t^2 \mu^4 (1+t k(y))^{-\frac94} \frac{3^\frac34 \sqrt{5}}{16}
\Big(\int_{\rand B_1(0)} \big|D^2k(y)(x)^2\big|^2\Big).
\end{align*}
\end{lemma}
\begin{proof}
We proceed as in Lemma \ref{l:expansion:3} and suppress the dependence of $w$ and $\vec{\alpha}$ 
on $t$, $\mu$, and $y$. From Lemmas \ref{l:estimates1} and \ref{p:implicit} we have for $0\le i \le 3$  
\begin{align}
\label{eq:23}
-5 t &\int \big(k(x)-T^1_{k,y}(x)\big)\,(z_{\mu,y}+w_0)^{4} (w-w_0-w_1) (\dot\xi_{\mu,y})_i \notag\\
&= -5 t \int \frac12 D^2 k(y)(x-y)^2 \,(z_{\mu,y}+w_0)^{4} w_2 (\dot\xi_{\mu,y})_i 
+ t^2 O_{A_0}(\mu^{4+\frac14}). 
\end{align}
Furthermore,
\begin{align}
\label{eq:29}
- 10 \int &(1+t k(x)) (z_{\mu,y}+w_0)^{3} (w-w_0-w_1)^2 (\dot \xi_{\mu,y})_i \notag\\
&= - 10 \int (1+t k(y)) (z_{\mu,y}+w_0)^{3} (w_2)^2 (\dot \xi_{\mu,y})_i +  t^2 O_{A_0}(\mu^{4+\frac14}).   
\end{align}
{\bf Case $1\le i \le 3$:}
Since $w_2$ is even and $(\dot \xi_{\mu,y})_i$ is odd for $1 \le i \le 3$ 
the integrals in \eqref{eq:23}-\eqref{eq:29} vanish.
Thus, from \eqref{eq:114} and the definition of $\vec{\a}_3$ we get
\begin{align*}
\Big(\vec{\alpha}-\sum_{j=1}^3\vec{\alpha}_j\Big)_i
&= -\frac{t\mu^4}{(1+t k(y))^{\frac54}} \frac{3^\frac34 8}{\pi\sqrt{5}}
\int \frac{\big(k(x+y)- T_{k(\cdot+y),0}^3(x)\big) x_i}{(\mu^2 +|x|^2)^{4}} \\
&\quad + t^2 O_{A_0}\big(|\nabla k(y)|^2\mu^{2}+|\nabla k(y)|\mu^{3}+\mu^{4+\frac14}\big).
\end{align*}
Since $D^4 k(y)(x)^4$ is even we may proceed analogously as in \eqref{eq:32} and prove
the claim of the lemma if $1\le i \le 3$.\\
{\bf Case $i=0$:}
From \eqref{eq:114}, the definition of $\vec{\a}_3$, \eqref{eq:74}, \eqref{eq:23}-\eqref{eq:29}, 
and the fact that $D^3 k(y)(x)^3$ is odd, 
we arrive at
\begin{align}
(\vec{\alpha}&-\vec{\alpha}_1-\vec{\alpha}_2-\vec{\alpha}_3)_0 \notag \\ 
&=\frac{3^\frac34 4t}{\pi\sqrt{5}(1+t k(y))^{\frac54}} \int \big(k(x+y)-T^3_{k(\cdot+y),0}(x)\big) \notag \\
&\qquad \qquad \qquad \qquad\qquad
\Bigg(\frac{3\big|\frac{x}{\mu}\big|^{-2}+3\big|\frac{x}{\mu}\big|^{-4}
+\big|\frac{x}{\mu}\big|^{-6}}{\mu^{3}\big(1+\big|\frac{x}{\mu}\big|^2\big)^{3}} 
+\frac{2}{\mu^{3}\big(1+\big|\frac{x}{\mu}\big|^2\big)^{4}}\Bigg)\notag \\
&\quad -5 t \int \frac{1}{2}D^2 k(y)(x-y)^2 \,(z_{\mu,y}+w_0)^{4} w_2 (\dot\xi_{\mu,y})_0\notag \\
&\quad -10 \int (1+t k(y)) (z_{\mu,y}+w_0)^{3}(w_2)^2(\dot\xi_{\mu,y})_0\notag\\
&\quad +t^2 O_{A_0}\big(|\nabla k(y)|^2\mu^{2}+|\nabla k(y)|\mu^{3}+\mu^{4+\frac14}\big).
\label{eq:73}
\end{align}
From Lemma \ref{l:estimates1}, \eqref{eq:A34}, and \eqref{eq:A35} we get
\begin{align*}
\int &\big(k(x+y)-T^3_{k(\cdot+y),0}(x)\big)\,
\Bigg(\frac{3\big|\frac{x}{\mu}\big|^{-2}+3\big|\frac{x}{\mu}\big|^{-4}
+\big|\frac{x}{\mu}\big|^{-6}}{\mu^{3}\big(1+\big|\frac{x}{\mu}\big|^2\big)^{3}} 
+\frac{2}{\mu^{3}\big(1+\big|\frac{x}{\mu}\big|^2\big)^{4}}\Bigg)\\
&= \mu^4 \int \frac{1}{4!}D^4 k(y)(x)^4 
\frac{5+6|x|^{-2}+4|x|^{-4}+|x|^{-6}}{(1+|x|^2)^{4}} 
+O_{A_0}\big(\mu^{4+\frac12}\big)\\
&= \mu^4 \frac{\pi^2}{24} \laplace^2 k(y) + O_{A_0}\big(\mu^{4+\frac12}\big).  
\end{align*}
For the second term in \eqref{eq:73} we obtain from Lemmas \ref{p:implicit} and \ref{l:comp:w2} 
\begin{align*}
-5 t &\int \frac{1}{2}D^2 k(y)(x-y)^2 \,(z_{\mu,y}+w_0)^{4} w_2 (\dot\xi_{\mu,y})_0\\
&= -  \frac{t^2 \mu^4 4 \sqrt{15}}{\pi(1+t k(y))^{\frac94}} 
\int \frac{1}{2}D^2 k(y)(x)^2 {\tilde w}_2(y)\, (1+|x|^2)^{-\frac52}\Big(1-\frac{2}{1+|x|^2}\Big).
\end{align*}
Moreover, from  Lemma \ref{l:comp:w2}, \eqref{eq:59}, and \eqref{eq:A43}
\begin{align}
\label{eq:85}
\int \frac{1}{2}&D^2 k(y)(x)^2 {\tilde w}_{2}(y) (1+|x|^2)^{-\frac52}\big(1-2(1+|x|^2)^{-1}\big) \notag \\
&= O_{A_0}\big(|\laplace k(y)|^2\big) 
+ \Big(\sum_{1\le l<m\le 3}\psi_{l,m}(y)^2 +\sum_{l=2}^3\psi_{l,l}(y)^2 \Big) \frac{3^\frac54 \sqrt{\pi}}{4}\notag\\ 
&\quad \times \sum_{j=0}^\infty \frac{(3+j)j!}{\Gamma(\frac72+j)} 
\int_0^\infty \frac{r^6}{(1+r^2)^{5}} \Big(1-\frac{2}{1+r^2}\Big)
{\mathcal P}_j^{(\frac52,\frac52)}\Big(1-\frac{2}{1+r^2}\Big)\notag\\
&= O_{A_0}\big(|\laplace k(y)|^2\big)
+ \frac{3^\frac54 \pi}{64} \Big(\int_{\rand B_1(0)} \big|D^2k(y)(x)^2-\frac{\laplace k(y)}{3}|x|^2\big|^2\Big) \notag\\
&\quad \times \sum_{j=0}^\infty \frac{2+6j+j^2}{(5+j)(4+j)(2+j)(1+j)} \notag\\
&= O_{A_0}\big(|\laplace k(y)|\big)
+ \frac{3^\frac54 5 \pi}{1024} \Big(\int_{\rand B_1(0)} \big|D^2k(y)(x)^2\big|^2\Big).  
\end{align}
For the remaining term in \eqref{eq:73} we derive
\begin{align*}
-10 \int &(1+t k(y)) (z_{\mu,y}+w_0)^{3}(w_2)^2(\dot\xi_{\mu,y})_0\\
&= -10 \frac{t^2 \mu^4}{(1+t k(y))^{\frac94}} \frac{3^{\frac34} 4}{\pi \sqrt{15}}
\int (1+|x|^2)^{-2}\Big(1-\frac{2}{1+|x|^2}\Big) ({\tilde w}_{2}(y))^2.
\end{align*}
By Lemma \ref{l:comp:w2}, \eqref{eq:85}, and ßß we have with  $\beta_j :=  \frac{(3+j)j!}{\Gamma(j+\frac{7}2)}$
\begin{align}
\label{eq:88}
\int (1&+|x|^2)^{-2}\Big(1-\frac{2}{1+|x|^2}\Big) ({\tilde w}_{2}(y))^2 \notag \\
&= O_{A_0}\big(|\laplace k(y)|\big) +
\frac14\Big(\int_{\rand B_1(0)} \big|D^2k(y)(x)^2\big|^2\Big) \frac{3^{\frac52} \pi}{16} \notag \\
&\qquad\times \int_0^\infty r^6 (1+r^2)^{-7} \Big(1-\frac{2}{1+r^2}\Big) 
\bigg(\sum_{j=0}^\infty \beta_j {\mathcal P}_j^{(\frac52,\frac52)}\Big(1-\frac{2}{1+r^2}\Big)\bigg)^2\notag \\
&= \frac14\Big(\int_{\rand B_1(0)} \big|D^2k(y)(x)^2\big|^2\Big) \frac{3^{\frac52} \pi}{16} \frac{1}{288}
+O_{A_0}\big(|\laplace k(y)|\big).  
\end{align}
Combining the computations in \eqref{eq:73}-~\eqref{eq:88} ends the proof.
\end{proof}

\begin{lemma}
\label{l:expansion:5}
Under the assumptions of Lemma \ref{p:implicit} suppose $k\in C^6(\rz^3)$ and 
$\|D^6 k\|_{\infty} \le A_0$. Then there holds
\begin{align*}
\|(\vec{\alpha}(&t,\mu,y))_0-\sum_{j=1}^5(\vec{\alpha}_{j}(t,\mu,y))_0\| 
= O_{A_0}(t \mu^{6})\\
&+ t^2 O_{A_0}\Big(\mu^2|\nabla k(y)|^2+\mu^3|\nabla k(y)|+ \mu^4 |\laplace k(y)|+ \mu^{4+\frac14}\Big).
\end{align*}
where $(\vec{\alpha}_{5}(t,\mu,y))_0$ is given by
\begin{align*}
(\vec{\alpha}_{5}(t,\mu,y))_0:=  -\frac{t \mu^5}{(1+t k(y))^{\frac54}} \frac{3^\frac34 4 \sqrt{5}}{\pi}
\cpvint{}{}  \Big(k(y+x)-T^4_{k(\cdot+y),0}(x) \Big) \frac1{|x|^8}.   
\end{align*}
\end{lemma}
\begin{proof}
From the proof of Lemma \ref{l:expansion:4}
we infer
\begin{align*}
(\vec{\alpha}(&t,\mu,y))_0-\sum_{j=1}^4(\vec{\alpha}_{j}(t,\mu,y))_0 
= -t (1+t k(y))^{-\frac54} \frac{3^\frac54 4}{\pi \sqrt{15}}\\
&\int \Big(k(y+x)-\sum_{\ell=0}^4 \frac{1}{\ell!} D^\ell k(y)(x)^\ell \Big) 
\frac{5+6\big|\frac{x}{\mu}\big|^{-2}+ 4\big|\frac{x}{\mu}\big|^{-4}+\big|\frac{x}{\mu}\big|^{-6}}
{\mu^3 \Big(1+\big|\frac{x}{\mu}\big|^{2}\Big)^4}\\
&+ t^2 O_{A_0}\big(\mu^2|\nabla k(y)|^2+\mu^3|\nabla k(y)|+ \mu^4 |\laplace k(y)|+ \mu^{4+\frac14}\big)
\end{align*}
Since $D^5k(y)(x)^5$ is odd, analogously as in Lemma \ref{l:expansion:3} we find
\begin{align*}
\int \Big(&k(y+x)-\sum_{\ell=0}^4 \frac{1}{\ell!} D^\ell k(y)(x)^\ell \Big) 
\frac{5+6\big|\frac{x}{\mu}\big|^{-2}+ 4\big|\frac{x}{\mu}\big|^{-4}+\big|\frac{x}{\mu}\big|^{-6}}
{\mu^3 \Big(1+\big|\frac{x}{\mu}\big|^{2}\Big)^4}\\
&= 5\mu^5 \cpvint{}{}  \Big(k(y+x)-\sum_{\ell=0}^4 \frac{1}{\ell!} D^\ell k(y)(x)^\ell \Big) \frac1{|x|^8}
+O_{A_0}(\mu^{6}),
\end{align*}
which ends the proof.
\end{proof}

\section{Derivatives of $\vec{\a}$}
\label{s:derivatives_alpha}

\begin{lemma}
\label{l:2nd_derivative_alpha}
Under the assumptions of Lemma \ref{p:implicit} we have for all
$(t,\mu,y) \in \Omega$ with $|\mu|\le 1$ and $1 \le i,j \le 3$
\begin{align}
\label{eq:48}
\Big|\frac{1}{t \mu}\frac{\rand \alpha(t,\mu,y)_i}{\rand y_j} 
&+ \frac{\pi}{3^\frac14 \sqrt{5}}(1+t k(y))^{-\frac54} \frac{\rand^2 k(y)}{\rand x_i \rand x_j}\Big| 
\le O_{A_0}\big(|\nabla k(y)|^2+\mu^{1+\frac14}\big),
\end{align}
\begin{align}
\label{eq:89}
\Big|\frac{1}{t \mu^2}\frac{\rand \alpha(t,\mu,y)_0}{\rand y_j} 
&+ \frac{\pi}{3^\frac14 \sqrt{5}} (1+t k(y))^{-\frac54}\frac{\rand}{\rand x_j} \laplace k(y)\Big| \notag \\
&\le O_{A_0}\big(|\nabla k(y)|^2\mu^{-1}+ \mu^{\frac14}\big).  
\end{align}
\end{lemma}
\begin{proof}
In the sequel we will suppress the dependence of $w$ and $\vec{\a}$ on $(t,\mu,y)$.  
Since $H(t,\mu,y,w,\vec{\a})\equiv 0$ we have
\begin{align}
\label{eq:28}
-\frac{\rand H}{\rand y} =\frac{\rand H}{\rand (w,\vec{\a})} 
\binom{\frac{\rand w}{\rand y}}{\frac{\rand \vec{\a}}{\rand y}}=
\frac{\rand H}{\rand w}\frac{\rand w}{\rand y}+ \frac{\rand H}{\rand \vec{\a} }\frac{\rand \vec{\a}}{\rand y}, 
\end{align}
where
\begin{align}
\label{eq:49}
\frac{\rand H}{\rand y_j}(t,\mu,y,w,\vec{\a}) = 
\Big(f_t''(z_{\mu,y}+w)\frac{\rand z_{\mu,y} }{\rand y_j} 
- \vec{\a} \cdot \frac{\rand \dot{\xi}_{\mu,y}}{\rand y_j}, 
\big(\langle w, \frac{\rand (\dot{\xi}_{\mu,y})_l}{\rand y_j} \rangle\big)_l  
\Big).  
\end{align}
A direct calculation gives for $0\le l\le 3$
\begin{align}
\label{eq:24}
\big\|\frac{\rand (\dot{\xi}_{\mu,y})_l}{\rand y_j} \big\|
+\big\|\frac{\rand z_{\mu,y}}{\rand y_j} \big\| \le {\rm const}\ \mu^{-1}.   
\end{align}
Differentiating $\langle z_{\mu,y}, (\dot{\xi}_{\mu,y})_l \rangle \equiv 0$ leads to
\begin{align*}
\langle z_{\mu,y}, \frac{\rand (\dot{\xi}_{\mu,y})_l}{\rand y_j} \rangle =  
\langle \frac{\rand z_{\mu,y} }{\rand y_j}, (\dot{\xi}_{\mu,y})_l \rangle,  
\end{align*}
and with \eqref{eq:14} and \eqref{eq:24} we arrive at
\begin{align}
\label{eq:25}
\Big|\langle w,\frac{\rand (\dot{\xi}_{\mu,y})_l}{\rand y_j}\rangle -
\big((1+t k(y))^{-\frac14} -1\big)\langle \frac{\rand z_{\mu,y} }{\rand y_j}, (\dot{\xi}_{\mu,y})_l \rangle \Big|
= tO_{A_0}(|\nabla k(y)|+\mu). 
\end{align}
By \eqref{eq:24}, the expansion of $\vec{\a}$ in Lemma \ref{p:implicit}, and (\ref{eq:77}) we see
\begin{align}
\label{eq:26}
\Big\|f_t''(z_{\mu,y}+w)\frac{\rand z_{\mu,y} }{\rand y_j}\Big\| +
\Big\| \vec{\a} \cdot \frac{\rand \dot{\xi}_{\mu,y}}{\rand y_j} \Big\| 
\le  tO_{A_0}(|\nabla k(y)|+\mu). 
\end{align}
From \eqref{eq:25}-\eqref{eq:26} we get 
\begin{align*}
\Big\|\frac{\rand H}{\rand(w,\vec\alpha)}(t,\mu,y,w,\vec{\a})&
\binom{\big((1+t k(y))^{-\frac14} -1\big)\frac{\rand z_{\mu,y} }{\rand y_j}}{0} + 
\frac{\rand H}{\rand y_j}(t,\mu,y,w,\vec{\a})\Big\| \\
&\qquad \qquad \qquad \le t O_{A_0}(|\nabla k(y)|+\mu)  
\end{align*}
which implies due to the uniform bound of the inverse (see Remark \ref{r:global})
\begin{align}
\label{eq:38}
\Big\|\frac{\rand w}{\rand y_j}-\big((1+t k(y))^{-\frac14} -1\big)\frac{\rand z_{\mu,y}}{\rand y_j} \Big\|
+ \Big\| \frac{\rand \vec{\a}}{\rand y_j} \Big\| \le tO_{A_0}(|\nabla k(y)|+\mu).
\end{align}
From \eqref{eq:91} and \eqref{eq:28}-~\eqref{eq:49} we deduce after testing with $(\dot{\xi}_{\mu,y})_j$
\begin{align}
\label{eq:51}
\frac{\rand \vec{\a}_j}{\rand y_i} &= 
f_t''(z_{\mu,y}+w)\big(\frac{\rand z_{\mu,y} }{\rand y_i}+ \frac{\rand w}{\rand y_i}\big)(\dot{\xi}_{\mu,y})_j
-\sum_{l=0}^3 \vec{\a}(t,\mu,y)_l \langle \frac{\rand (\dot{\xi}_{\mu,y})_l}{\rand y_i},(\dot{\xi}_{\mu,y})_j \rangle. 
\end{align}
From Lemma \ref{l:estimates1}, \eqref{eq:77}, \eqref{eq:38}, and the fact that 
$(\dot{\xi}_{\mu,y})_j \in N(f_0''(z_{\mu,y}))$ we obtain
\begin{align}
\label{eq:53}
f_t''(z_{\mu,y}&+w) \big(\frac{\rand z_{\mu,y} }{\rand y_i}+ \frac{\rand w}{\rand y_i}\big)(\dot{\xi}_{\mu,y})_j \notag\\
&= -\frac{5 t}{(1+t k(y))^{\frac54}}\int \Big(\sum_{\ell=1}^3\frac{1}{\ell!}D^\ell k(y)(x-y)^\ell\Big) (z_{\mu,y})^4 
\big(\frac{\rand z_{\mu,y} }{\rand y_i} \big) (\dot{\xi}_{\mu,y})_j\notag\\
&\quad -20 \int (z_{\mu,y})^3 (w-w_0) 
\big(\frac{\rand z_{\mu,y} }{\rand y_i} \big) (\dot{\xi}_{\mu,y})_j \notag \\
&\quad+O_{A_0}\big(t\mu^{2+\frac34}+t^2(|\nabla k(y)|^2\mu+|\nabla k(y)|\mu^\frac74+\mu^{2+\frac34})\big)
\end{align}
Differentiating the identity $f_0''(z_{\mu,y})(\dot{\xi}_{\mu,y})_j=0$ with respect to $y_i$
leads to
\begin{align*}
0= f_0'''(z_{\mu,y})\big(\frac{\rand z_{\mu,y} }{\rand y_i} \big)(\dot{\xi}_{\mu,y})_j
+ f_0''(z_{\mu,y})\big(\frac{\rand (\dot{\xi}_{\mu,y})_j }{\rand y_i} \big), 
\end{align*}
and we get from Lemma \ref{l:estimates1}, \eqref{eq:13}, and \eqref{eq:72}
\begin{align*}
-&20 \int (z_{\mu,y})^3 (w-w_0)\frac{\rand z_{\mu,y} }{\rand y_i} (\dot{\xi}_{\mu,y})_j \notag \\
&=  f_0'''(z_{\mu,y})(w-w_0) \frac{\rand z_{\mu,y} }{\rand y_i} (\dot{\xi}_{\mu,y})_j 
= -f_0''(z_{\mu,y}) \frac{\rand (\dot{\xi}_{\mu,y})_j }{\rand y_i}(w-w_0) \notag\\
&= -\sum_{l=0}^3 \big(\vec{\a}(t,\mu,y)\big)_l
\big\langle (\dot{\xi}_{\mu,y})_l,\frac{\rand (\dot{\xi}_{\mu,y})_j}{\rand y_i}\big\rangle\\
&\quad -\frac{t}{(1+t k(y))^{\frac54}} \int \sum_{\ell=1}^3 \frac{1}{\ell!} D^\ell k(y)(x-y)^\ell (z_{\mu,y})^5
\frac{\rand (\dot{\xi}_{\mu,y})_j}{\rand y_i} \notag \\
&\quad + t O_{A_0}\big(\mu |\nabla k(y)|^2+\mu^{2+\frac14}\big).
\end{align*}
Differentiating $\langle (\dot{\xi}_{\mu,y})_l,(\dot{\xi}_{\mu,y})_j\rangle \equiv const$ with respect
to $y_i$ we obtain
\begin{align*}
\big\langle \frac{\rand (\dot{\xi}_{\mu,y})_l}{\rand y_i},(\dot{\xi}_{\mu,y})_j,\big\rangle 
= - \big\langle (\dot{\xi}_{\mu,y})_l,\frac{\rand (\dot{\xi}_{\mu,y})_j}{\rand y_i}\big\rangle.
\end{align*}
Inserting the above computations in \eqref{eq:53} and \eqref{eq:51} leads to
\begin{align*}
\frac{\rand \vec{\a}_j}{\rand y_i} 
&=  -\frac{t}{(1+t k(y))^{\frac54}}\int \Big(\sum_{\ell=1}^3\frac{1}{\ell!}D^\ell k(y)(x-y)^\ell\Big) 
\frac{d}{d y_i}\Big((z_{\mu,y})^5 (\dot{\xi}_{\mu,y})_j\Big)\\
&\quad + tO_{A_0}\big(|\nabla k(y)|^2\mu+\mu^{2+\frac14}\big).
\end{align*}
As $\frac{\rand}{\rand y_i}z_{\mu,y}(x) = -\frac{\rand}{\rand x_i}z_{\mu,y}(x)$ and 
$\frac{\rand}{\rand y_i}(\dot{\xi}_{\mu,y})_j(x) = -\frac{\rand}{\rand x_i}(\dot{\xi}_{\mu,y})_j(x)$
we obtain by partial integration
\begin{align*}
\int \Big(\sum_{\ell=1}^3&\frac{1}{\ell!}D^\ell k(y)(x-y)^\ell\Big) 
\frac{d}{d y_i}\Big((z_{\mu,y})^5 (\dot{\xi}_{\mu,y})_j\Big)\\ 
&= \int \frac{d}{d x_i}\Big(\sum_{\ell=1}^3\frac{1}{\ell!}D^\ell k(y)(x-y)^\ell\Big) 
(z_{\mu,y})^5 (\dot{\xi}_{\mu,y})_j\\
&= \int \Big(\sum_{\ell=0}^2\frac{1}{\ell!}D^\ell \frac{\rand k}{\rand x_i}(y)(x-y)^\ell\Big) 
(z_{\mu,y})^5 (\dot{\xi}_{\mu,y})_j.
\end{align*}
The latter integral may be evaluated as in Lemma \ref{p:implicit} and yields the claim.
\end{proof} 

\begin{lemma}
\label{l:diff_alpha_0_eps}
Under the assumptions of Lemma \ref{l:2nd_derivative_alpha} we have
\begin{align}
\label{eq:86}
\Big\|\frac{\rand w}{\rand t}- \sum_{i=0}^2\frac{\rand w_i}{\rand t}\Big\|
+ \Big\|\frac{\rand \vec{\a}}{\rand t}- \sum_{i=1}^2\frac{\rand \vec{\a_i}}{\rand t}\Big\|
=  O_{A_0}\big(|\nabla k(y)|^2 \mu^2+\mu^{2+\frac14}\big).
\end{align}  
\end{lemma}
\begin{proof}
In the sequel we will suppress the dependence of $w$ and $\vec{\a}$ on $(t,\mu,y)$.  
Since $H(t,\mu,y,w,\vec{\a})\equiv 0$ we have
\begin{align}
\label{eq:83}
\Big(\int k(x) &(z_{\mu,y}+w)^5 \cdot, \vec{0}\Big)
= -\frac{\rand H}{\rand t}\bigg|_{(t,\mu,y,w,\vec{\a})} 
= \frac{\rand H}{\rand (w,\vec{\a})}\bigg|_{(t,\mu,y,w,\vec{\a})}
\binom{\frac{\rand w}{\rand t}}{\frac{\rand \vec{\a}}{\rand t}}\notag \\
&=\Big(f_t''(z_{\mu,y}+w)\frac{\rand w}{\rand t} 
- \frac{\rand \vec{\a}}{\rand t} \cdot \dot{\xi}_{\mu,y}, 
\big(\langle \frac{\rand w}{\rand t},(\dot{\xi}_{\mu,y})_l\rangle\big)_l 
\Big).
\end{align} 
Differentiating the identities
\begin{align*}
f_0'(z_{\mu,y}+w_0)&= t \int k(y)(z_{\mu,y}+w_0)^5 \cdot,\\
f_0''(z_{\mu,y})\sumg_{i=1}^2w_i -\sumg_{i=1}^2 \vec{\a}_i \cdot \dot{\xi}_{\mu,y} &=
t \int\Big(\sum_{\ell=1}^2 \frac{1}{\ell!}D^\ell k(y)(x-y)^\ell\Big)(z_{\mu,y}+w_0)^5 \cdot. 
\end{align*}
with respect to $t$ leads to
\begin{align}
\label{eq:78}
&f_0''(z_{\mu,y})\sum_{i=0}^2\frac{\rand  w_i}{\rand t} 
-\sum_{i=1}^2\frac{\rand \vec{\a}_i}{\rand t} \cdot \dot{\xi}_{\mu,y}
= \int\Big(\sum_{l=0}^2 \frac{1}{\ell!}D^\ell k(y)(x-y)^\ell\Big)(z_{\mu,y}+w_0)^5 \cdot \notag \\
&\qquad \quad 
+5t \int\Big(\sum_{\ell=1}^2 \frac{1}{\ell!}D^\ell k(y)(x-y)^\ell\Big)(z_{\mu,y}+w_0)^4\frac{\rand w_0}{\rand t} \cdot.
\end{align}
Furthermore, we note that
\begin{align}
\label{eq:87}
\frac{\rand  w_0}{\rand t} = -\frac{k(y)}{4}(1+t k(y))^{-1}(z_{\mu,y}+w_0).  
\end{align}
For $\frac{\rand  w_i}{\rand t}= O_{A_0}(\mu^i)$ as $\mu \to 0$ we get
from Lemma \ref{l:estimates1}, \eqref{eq:77}, \eqref{eq:78}, and \eqref{eq:87}
\begin{align}
\label{eq:79}
f_t''&(z_{\mu,y}+w)\sum_{i=0}^2\frac{\rand  w_i}{\rand t} 
-\frac{\rand (\vec{\a}_1+\vec{\a}_2)}{\rand t} \cdot \dot{\xi}_{\mu,y} \notag \\
&= \int\Big(\sum_{l=0}^2 \frac{1}{\ell!}D^\ell k(y)(x-y)^\ell\Big)(z_{\mu,y}+w_0)^5 \cdot \notag \\
&\quad + 5 \int k(y)(z_{\mu,y}+w_0)^4(w-w_0)\cdot \ +t O_{A_0}\big(|\nabla k(y)|^2 \mu^2+\mu^\frac94\big).
\end{align}
Moreover, by Lemmas \ref{l:estimates1} and \ref{p:implicit}
\begin{align}
\int &k(x) (z_{\mu,y}+w)^5 \cdot \notag\\ 
&= \int\Big(\sum_{l=0}^2 \frac{1}{\ell!}D^\ell k(y)(x-y)^\ell\Big)(z_{\mu,y}+w_0)^5 \cdot 
+ O_{A_0}(\mu^\frac94)\notag \\
\label{eq:80}
&\quad +5 \int k(y) (z_{\mu,y}+w_0)^4(w-w_0) \cdot
+t O_{A_0}\big(|\nabla k(y)|^2 \mu^2+\mu^3\big).
\end{align}
Combining \eqref{eq:79}, \eqref{eq:80} and the fact that $\frac{\rand w_i}{\rand t}$ remains 
in $T_{z_{\mu,y}}Z^\perp$ we get
\begin{align*}
\frac{\rand H}{\rand (w,\vec{\a})}
\binom{\sum_{i=0}^2\frac{\rand w_i}{\rand t}}{\sum_{i=1}^2\frac{\rand \vec{\a_i}}{\rand t}}
&= \frac{\rand H}{\rand t}
+ O_{A_0}\big(|\nabla k(y)|^2 \mu^2+ \mu^{2+\frac14}\big).
\end{align*}
and the claim of the lemma follows from Remark \ref{r:global}. 
\end{proof}

\begin{lemma}
\label{l:d_alpha_0_d_eps_4}
Under the assumptions of Lemma \ref{l:2nd_derivative_alpha} we have
\begin{align}
\label{eq:98}
\frac{\rand (\vec{\a})_0}{\rand t}
&= \frac1t(\vec{\a})_0+ \sum_{j=2}^4\frac{\rand (\vec{\a_j})_0}{\rand t}-\frac1t(\vec{\a_j})_0 \notag \\
&\quad + tO_{A_0}\big(|\nabla k(y)|^2 \mu^2+ |\nabla k(y)|\mu^3+\mu^4 |\laplace k(y)|^2+ \mu^{4+\frac14}\big),
\end{align}
and for $1 \le i \le 3$
\begin{align}
\label{eq:99}
\frac{\rand (\vec{\a})_i}{\rand t}
&= \frac1t(\vec{\a})_i+ \sum_{j=1}^3\frac{\rand (\vec{\a_j})_i}{\rand t}-\frac1t(\vec{\a_j})_i \notag \\
&\quad + tO_{A_0}\big(|\nabla k(y)|^2 \mu^2+ |\nabla k(y)|\mu^3+ \mu^{3+\frac12}\big).
\end{align}
\end{lemma}
\begin{proof}
By \eqref{eq:83} we have for $0 \le i \le 3$ 
\begin{align*}
\frac{\rand (\vec{\a})_i}{\rand t} = f_t''(z_{\mu,y}+w) \frac{\rand w}{\rand t} (\dot{\xi}_{\mu,y})_i
- \int k(x)(z_{\mu,y}+w)^5 (\dot{\xi}_{\mu,y})_i.
\end{align*}
To prove the claim of the lemma we will proceed termwise.
In the calculations below certain terms will vanish simply because we
are integrating a product of an odd and an even function. Moreover, we often use Lemma \ref{l:estimates1}
without mentioning it explicitly.\\
For $(\dot{\xi}_{\mu,y})_i \in N(f_0''(z_{\mu,y})$ and by \eqref{eq:77} we see
\begin{align*}
f_t''(z_{\mu,y}+w) \frac{\rand w}{\rand t} (\dot{\xi}_{\mu,y})_i &=
-5 t \int(k(x)-k(y))(z_{\mu,y}+w_0)^4 \frac{\rand w}{\rand t} (\dot{\xi}_{\mu,y})_i\\
&\quad -20 \int(1+t k(x))(z_{\mu,y}+w_0)^3(w-w_0) \frac{\rand w}{\rand t} (\dot{\xi}_{\mu,y})_i\\
&\quad -30 \int(1+t k(x))(z_{\mu,y}+w_0)^2(w-w_0)^2 \frac{\rand w}{\rand t} (\dot{\xi}_{\mu,y})_i\\
&\quad + t^3  O_{A_0}\big(|\nabla k(y)|^3 \mu^3+\mu^6\big).  
\end{align*}
Due to Lemma \ref{l:diff_alpha_0_eps} we may replace $\frac{\rand w}{\rand t}$ by 
$\sum_{i=0}^2 \frac{\rand w_i}{\rand t}$.  
By \eqref{eq:87} we obtain
\begin{align*}
-5 t &\int(k(x)-k(y))(z_{\mu,y}+w_0)^4 \frac{\rand w}{\rand t} (\dot{\xi}_{\mu,y})_i \notag\\
&= \frac{5 k(y)}{4(1+t k(y))}t \int(k(x)-k(y))(z_{\mu,y}+w_0)^5(\dot{\xi}_{\mu,y})_i \notag\\
&\quad -5t \int \frac12 D^2k(y)(x-y)^2 (z_{\mu,y}+w_0)^4 \frac{\rand w_2}{\rand t} (\dot{\xi}_{\mu,y})_i\notag \\
&\quad +t O_{A_0}\big(\mu^2 |\nabla k(y)|^2+\mu^3 |\nabla k(y)| +\mu^{4+\frac14}\big).
\end{align*}
From \eqref{eq:87} and as $(\dot{\xi}_{\mu,y})_i \in N(f_0''(z_{\mu,y}))$ 
and $w-w_0 \in T_{z_{\mu,y}}T^\perp$  
we see
\begin{align*}
-20 \int(1+t k(y))(z_{\mu,y}+w_0)^3&(w-w_0) \frac{\rand w_0}{\rand t} (\dot{\xi}_{\mu,y})_i\\
&=  \frac{k(y) \langle (w-w_0),(\dot{\xi}_{\mu,y})_i \rangle}{1+t k(y)} =0.  
\end{align*}
By Lemma \ref{l:estimates1} and \eqref{eq:86} we now get
\begin{align*}
-20 &\int(1+t k(x))(z_{\mu,y}+w_0)^3(w-w_0) \frac{\rand w}{\rand t} (\dot{\xi}_{\mu,y})_i \notag \\  
&= -20 \int(1+t k(y))(z_{\mu,y}+w_0)^3w_2 \frac{\rand w_2}{\rand t} (\dot{\xi}_{\mu,y})_0\notag \\
&\quad -20 t \int\frac12 D^2k(y)(x-y)^2 (z_{\mu,y}+w_0)^3 w_2 
\frac{\rand w_0}{\rand t} (\dot{\xi}_{\mu,y})_i \notag\\
&\quad + tO_{A_0}\big(\mu^2 |\nabla k(y)|^2+\mu^3 |\nabla k(y)|+\mu^{4+\frac14}\big).
\end{align*}
Furthermore, we see
\begin{align*}
-30 \int&(1+t k(x))(z_{\mu,y}+w_0)^2(w-w_0)^2 \frac{\rand w}{\rand t} (\dot{\xi}_{\mu,y})_i \notag \\ 
&= -30 \int (1+t k(y))(z_{\mu,y}+w_0)^2 (w_2)^2 \frac{\rand w_0}{\rand t} (\dot{\xi}_{\mu,y})_i\notag \\
&\quad +t^2 O_{A_0}\big(\mu^2 |\nabla k(y)|^2+\mu^3|\nabla k(y)| +\mu^{4+\frac14}\big)
\end{align*}
Since $z_{\mu,y}$ and $w$ are orthogonal to $(\dot{\xi}_{\mu,y})_i$ and 
$(\dot{\xi}_{\mu,y})_i \in N(f_0''(z_{\mu,y}))$, we may estimate using Lemma \ref{p:implicit}
\begin{align*}
-\int k(x) (z_{\mu,y}+w)^5(\dot{\xi}_{\mu,y})_i 
&= \frac{1}{t} f_t'(z_{\mu,y}+w)(\dot{\xi}_{\mu,y})_i 
+\frac{1}{t} \int (z_{\mu,y}+w)^5(\dot{\xi}_{\mu,y})_i \notag \\
&= \frac{1}{t}(\vec{\a})_i + \frac{10}{t} \int (z_{\mu,y}+w_0)^3(w_2)^2 (\dot{\xi}_{\mu,y})_i \notag \\
&\quad +t^2 O_{A_0}(|\nabla k(y)|^2 \mu^2+ |\nabla k(y)|\mu^3+\mu^{4+\frac14}).
\end{align*}
As
\begin{align*}
\frac{\rand w_2}{\rand t}= \Big(\frac1t -\frac{5 k(y)}{4(1+t k(y))}\Big)w_2,  
\end{align*}
we end up with integrals that are, up to a factor, computed in Section \ref{s:expansion}.
Summing up the results will give the claim of the lemma.
\end{proof}

\section{Solvability of $\vec{\a}(t,\mu,y)=0$}
\label{s:solv_alpha_equal_null}

\begin{lemma}
\label{l:implicit:curve}
Under the assumptions of Lemma \ref{p:implicit} suppose $y_0$ is a nondegenerate critical point of $k$, i.e.
\begin{align*}
\nabla k(y_0)=0 \text{ and } D^2k(y_0) \text{ is invertible, with }  \big\|\big(D^2k(y_0)\big)^{-1}\big\| \le A_0.
\end{align*}
Moreover, assume $\laplace k(y_0)=0$.
Consider the function $\hat{\alpha}$, defined by
\begin{align*}
\hat{\alpha}(t,\mu,y) := \frac{3^\frac14 \sqrt{5}}{t \mu \pi} 
(1+t k(y_0))^{\frac54}(\vec{\alpha}(t,\mu,y)_{1},\dots,\vec{\alpha}(t,\mu,y)_{3})^T,
\end{align*}
which is well defined and continuous in $\Omega$ (see Remark \ref{r:1_t_alpha_continuous}),
analogously we define $\hat{\alpha}_j(t,\mu,y)$.
Then there are $\delta_1= \delta_1(A_0)>0$ and a $C^2$-function $\beta$,
\begin{align*}
\beta: \{(t,\mu) \where t \in [-B_0,B_1],\, 0<\mu<\delta_1\} \to \rz^3,  
\end{align*}
such that 
\begin{align*}
\hat{\alpha}(t,\mu,\beta(t,\mu)) =0 \text{ for all } t \in [-B_0,B_1],\, 0<\mu<\delta_1,\\
\end{align*}
and
\begin{align*}
\beta(t,\mu) &= y_0
+\big(D^2k(y_0)\big)^{-1}\Big(\sum_{j=3}^4 \hat{\alpha}_j(t,\mu,y_0)\Big)
+ O_{A_0}(\mu^{3+\frac14}).  
\end{align*}
Moreover, $\beta$ is unique in the sense that, 
 if $y \in B_{\delta_1}(y_0)$ satisfies $\hat{\alpha}(t,\mu,y) =0$ 
for some $t \in [-B_0,B_1]$ and $0<\mu<\delta_1$, then $y= \beta(t,\mu)$. 
\end{lemma}
\begin{proof}
In view of \eqref{eq:48} we would like to apply the implicit function theorem to the
function $(\vec{\alpha}(t,\mu,y))_{1\le i \le 3}$ in the 
point $(t,0,y_0)$, but unfortunately $\vec{\alpha}$ may not be differentiable for $\mu=0$.
Instead we mimic the proof of the implicit function theorem and apply Banach's fixed-point 
theorem to the function 
\begin{align*}
F(t,\mu,y) := y + \big(D^2k(y_0)\big)^{-1}\hat{\alpha}(t,\mu,y)
\end{align*}
in $B_\delta(y_0)$, where $\delta>0$ will be chosen later. 
Fix $y_1,y_2 \in B_\delta(y_0)$, then by Lemma \ref{l:2nd_derivative_alpha}
\begin{align*}
|F(t,&\mu,y_1)-F(t,\mu,y_2)|\\
&= \big|(y_1-y_2)+ \big(D^2k(y_0)\big)^{-1}  
\int_0^1 \frac{\rand \hat{\alpha}}{\rand y}(t,\mu,y_2+t(y_1-y_2)) (y_1-y_2)\di t \big|\\
&\le \big|(y_1-y_2)- \bigg(\int_0^1 \big(D^2k(y_0)\big)^{-1} D^2k(y_2+t(y_1-y_2)) \di t\bigg) (y_1-y_2)\big|\\
&\qquad +  O_{A_0}\Big(\supm_{y \in B_\delta(y_0)}|\nabla k(y)|+ \mu^\frac14\Big) |y_1-y_2|\\
&\le O_{A_0}\Big(\delta + \mu^\frac14\Big) |y_1-y_2|.  
\end{align*}
For $y \in B_\delta(y_0)$ we estimate using Lemma \ref{p:implicit}
\begin{align*}
|F(t,\mu,y) -  y_0| &= \big|y-y_0 + \big(D^2k(y_0)\big)^{-1}\big(\hat{\alpha}(t,\mu,y)\big)\big|\\
&\le \Big|y-y_0 - \big(D^2k(y_0)\big)^{-1}\Big(\nabla k(y) + O_{A_0}(\mu^{2})\Big)\Big|\\
&\le O_{A_0}(\delta^2 +\mu^2).  
\end{align*}
Consequently, there is $\delta_1=\delta_1(A_0)>0$ such that 
$F(t,\mu,\cdot)$ is a contraction in $B_{\delta_1}(y_0)$ for any $0 <\mu < \delta_1$ and $t\in [-B_0,B_1]$.
From Banach's fixed-point theorem we may define $\beta(t,\mu)$ to be the unique fixed-point of $F(t,\mu,\cdot)$
in $B_{\delta_1}(y_0)$. After shrinking $\delta_1$ if necessary we may apply 
Lemma \ref{l:2nd_derivative_alpha} and the usual implicit function theorem
to see that the function $\beta$ is twice differentiable for $\mu>0$.\\ 
To deduce the expansion for small $\mu$ we fix $\rho>0$ and
\begin{align*}
y \in U_\rho:= B_\rho\bigg(y_0 +\big(D^2k(y_0)\big)^{-1}\Big(\sum_{j=3}^4 \hat{\alpha}_j(t,\mu,y_0)\Big)\bigg).
\end{align*}
Then, by Lemmas \ref{l:expansion:3} and \ref{l:expansion:4} 
\begin{align*}
\Big|F(t,&\mu,y) -  y_0 - \big(D^2k(y_0)\big)^{-1}\Big(\sum_{j=3}^4 \hat{\alpha}_j(t,\mu,y_0)\Big)\Big|\\
&\le \Big|y-y_0 - \big(D^2k(y_0)\big)^{-1}
\Big(
\nabla k(y)
+O_{A_0}(\mu\rho^2+\mu^2\rho+\mu^{3+\frac14})\Big)\Big|\\
&\le O_{A_0}\big(\rho^2 +\mu^2 \rho +\mu^{3+\frac14}\big).
\end{align*}
Hence, we may choose for small $\mu$ a radius $0<\rho= O_{A_0}(\mu^{3+\frac14})$ such that   
$F$ maps $U_\rho$ into itself and $U_\rho \subset B_{\delta_1}(y_0)$.
Consequently, the unique fixed-point $\beta(t,\mu)$ must lie in this ball. This ends the proof.
\end{proof}

\begin{lemma}
\label{l:expansion_tilde_alpha_0} 
Under the assumptions of Lemma \ref{l:implicit:curve}, if moreover $k \in C^6(\rz^3)$ and 
$\|D^6 k\|_\infty \le A_0$, then we have
\begin{align*}
\big(\a(t,\mu,\beta(t,\mu))\big)_0
&= -t \mu^3 (1+t k(y_0))^{-\frac54} \frac{3^\frac34 4}{\pi \sqrt{5}} a_0(y_0)\\
&\quad +t \mu^4 (1+t k(y_0))^{-\frac94} \frac{\pi 3^\frac34 \sqrt{5}}{30} \Big(a_1(y_0)+t a_2(y_0)\Big)\\
&\quad +t \mu^5 (1+t k(y_0))^{-\frac54} \frac{\pi 3^\frac34 \sqrt{5}}{30}  a_3(y_0)\\
&\quad + O_{A_0}(t\mu^{5+\frac12}+t^2\mu^{4+\frac14}),    
\end{align*}
where $a_i(y_0)=a_i(\theta)$ given in \eqref{eq:31} and \eqref{eq:33} with $k_\theta=k(\cdot+y_0)$.
If the assumption $k\in C^6(\rz^3)$ is dropped then the terms of order higher than $4$ in $\mu$
have to be replaced by $O_{A_0}(t\mu^{4+\frac14})$.  
\end{lemma}
\begin{proof}
In view of Lemma \ref{l:implicit:curve} and because $\nabla k(y_0)=0$
we may estimate functions of $y:= \beta(t,\mu)$ and of $k(y)=k(\beta(t,\mu))$ as follows
\begin{align}
\label{eq:105}
\begin{split}
&F(y)= F(y_0)+ F'(y_0) \big(D^2k(y_0)\big)^{-1}\Big(\sum_{j=3}^4 \hat{\alpha}_j(t,\mu,y_0)\Big)
+ O_{A_0}(\mu^{3+\frac14}),\\
&F(k(y))= F(k(y_0))+ O_{A_0}(\mu^{4}).     
\end{split}
\end{align}
To prove the claim of the lemma we expand $\a(t,\mu,\beta(t,\mu))_0$ according to Lemma \ref{l:expansion:5} 
and use \eqref{eq:105}. For instance we have
\begin{align*}
(\vec{\alpha}_2(&t,\mu,y))_0 
=  -t\mu^2 (1+t k(y_0))^{-\frac{5}{4}} \frac{\pi}{3^{\frac14}\sqrt{5}}\laplace k(y)+O_{A_0}(t\mu^6)\\
&=  t\mu^4 (1+t k(y_0))^{-\frac{5}{4}} \Bigg(\frac{\pi}{3^{\frac14}2\sqrt{5}} 
\nabla\laplace k(y_0) \Big(D^2k(y_0)\Big)^{-1}\nabla(\laplace k(y_0))\\
&\quad +\mu \frac{3^{\frac34}8}{\pi\sqrt{5}} 
\nabla\laplace k(y_0)\cdot \big(D^2k(y_0)\big)^{-1}
\cpvint{}{} \frac{\big(k(x+y_0)- T^3_{k(\cdot+y_0),0}(x)\big)x_i}{|x|^{8}}\Bigg)\\
&\quad +O_{A_0}(t\mu^{5+\frac14}).
\end{align*}
If we continue expanding the remaining terms given in Lemma \ref{l:expansion:5} the claim follows.
\end{proof}

\begin{lemma}
\label{l:d_gamma_d_eps}
Under the assumptions of Lemma \ref{l:implicit:curve} let $\laplace k(y_0)=0=a_0(y_0)$
and define 
\begin{align*}
\gamma(t,\mu):= \frac{1}{t\mu^{4}} (1+t k(y_0))^{\frac94} \frac{30}{\pi 3^\frac34 \sqrt{5}} 
\big(\a(t,\mu,\beta(t,\mu))\big)_0.
\end{align*}
Then
\begin{align*}
\frac{\rand \gamma(t,\mu)}{\rand t} &= a_2(y_0)+ O_{A_0}(\mu^{\frac14}).
\end{align*}
\end{lemma}
\begin{proof}
We have
\begin{align*}
\frac{d\big(\vec{\a}(t,\mu,\beta(t,\mu))\big)_0 }{d t} 
= \frac{\rand (\vec{\a})_0}{\rand t}\Big|_{(t,\mu,\beta(t,\mu))}
+ \frac{\rand (\vec{\a})_0}{\rand y}\Big|_{(t,\mu,\beta(t,\mu))} \frac{\rand \beta}{\rand t}\Big|_{(t,\mu)}. 
\end{align*}
The derivatives of $(\vec{\a})_0$ are computed in \eqref{eq:89} and \eqref{eq:98}. In order to compute 
the derivative of $\beta$ we use the fact that 
\begin{align*}
\tilde{\a}(t,\mu,\beta):= \big((\vec{\a}(t,\mu,\beta))_1,\dots ,(\vec{\a}(t,\mu,\beta))_3\big)^T = \vec{0}.   
\end{align*}
By \eqref{eq:105} and Lemmas  \ref{l:2nd_derivative_alpha} and \ref{l:d_alpha_0_d_eps_4} we have
\begin{align*}
\frac{\rand \beta}{\rand t}\Big|_{(t,\mu)} &= 
-\bigg(\frac{\rand \tilde{\a}}{\rand y}\Big|_{(t,\mu,\beta(t,\mu))}\bigg)^{-1}
\frac{\rand \tilde{\a}}{\rand t}\Big|_{(t,\mu,\beta(t,\mu))}\\
&=t^{-1} \mu^{-1} \frac{3^\frac14 \sqrt{5}}{\pi}(1+t k(y_0))^\frac54 \Big(\big(D^2k(y_0)\big)^{-1}+
O_{A_0}(\mu^{\frac14})\Big) \notag \\ 
&\quad \bigg[\frac1t \big(\vec{\a}(t,\mu,\beta)\big)_i- 
O_{A_0}\Big(\sum_{j=1}^3\big(\vec{\a}_j(t,\mu,\beta)\big)_i\Big)+ tO_{A_0}(\mu^{3+\frac12})\bigg]_{i=1\dots 3}\\
&= O_{A_0}(\mu^{2+\frac12}),
\end{align*}
where we used the fact that as $\vec{\a}(t,\mu,\beta)_i\equiv 0$ for $1\le i\le3$  
\begin{align*}
\sum_{j=1}^3 \big(\vec{\a}_j(t,\mu,\beta)\big)_i = O_{A_0}(t\mu^{3+\frac12}). 
\end{align*}
From \eqref{eq:89} we get
\begin{align*}
\frac{\rand (\vec{\a})_0}{\rand y}&\Big|_{(t,\mu,\beta(t,\mu))} \frac{\rand \beta}{\rand t}\Big|_{(t,\mu)}
= O_{A_0}(t\mu^{4+\frac12}).
\end{align*}
Furthermore, by Lemmas \ref{l:d_alpha_0_d_eps_4} 
\begin{align*}
\frac{d \vec{\a}(t,\mu,\beta)_0 }{dt}
&= \frac{1}{t} \vec{\a}(t,\mu,\beta)_0 + \sum_{j=2}^4\frac{\rand \vec{\a_j}(t,\mu,\beta)_0}{\rand t}
-\frac1t \vec{\a_j}(t,\mu,\beta)_0\\
&\quad +O_{A_0}(t\mu^{4+\frac14}).
\end{align*}
The definition of $\gamma$, \eqref{eq:105}, and Lemma \ref{l:expansion_tilde_alpha_0}  
yield the claim.
\end{proof}

\begin{lemma}
\label{l:implicit:t(mu)}
Under the assumptions of Lemma \ref{l:d_gamma_d_eps} suppose $a_2(y_0)\not= 0$ and
either $a_1(y_0) \neq 0$ or ($a_1(y_0)=0$ and $a_3(y_0) \neq 0$). Moreover let
\begin{align*}
A_0 \ge 
\begin{cases}
{|a_1(y_0)|^{-1}}+2 |a_1(y_0)||a_2(y_0)|^{-1}+|a_2(y_0)|^{-1} &\text{ if } a_1(y_0) \neq 0,\\
|a_2(y_0)|^{-1}+|a_3(y_0)|^{-1}+ |a_3(y_0)| &\text{ if } a_1(y_0) = 0, 
\end{cases}
\end{align*}
and $-\frac{a_1(y_0)}{a_2(y_0)} \in (-B_0,B_1)$.
Then there exist $\delta_2=\delta_2(A_0)>0$ and a $C^1$-function $\tilde{t}$,
\begin{align*}
\tilde{t}: \{\mu \where 0<\mu<\delta_2\} \to (-B_0,B_1) \setminus \{0\},  
\end{align*}
such that $\big(\a(\tilde{t}(\mu),\mu,\beta(\tilde{t}(\mu),\mu))\big)_0\equiv 0$ for all $0<\mu<\delta_2$
and 
\begin{align}
\label{eq:102}
\tilde{t}(\mu)=-\frac{1}{a_2(y_0)} 
\begin{cases}
a_1(y_0)+ O_{A_0}(\mu^\frac14) &\text{ if } a_1(y_0) \neq 0\\
a_3(y_0)\mu + O_{A_0}(\mu^{1+\frac14}) &\text{ if } a_1(y_0) = 0.
\end{cases}
\end{align}
Moreover $\tilde{t}$ is unique in the sense that, if $t \in (-B_0,B_1)$ and $0<\mu<\delta_2$ 
satisfy $\big(\a({t},\mu,\beta({t},\mu))\big)_0= 0$ then $t= \tilde{t}(\mu)$.
\end{lemma}
\begin{proof}
We only sketch the proof, which is similar to the proof of Lemma \ref{l:implicit:curve}. 
We will apply Banach's fixed-point theorem to the function
\begin{align*}
F_\mu(t)=F(t,\mu):= t- a_2(y_0)^{-1}\gamma(t,\mu),  
\end{align*}
where $\gamma$ is given in Lemma \ref{l:d_gamma_d_eps}. 
To this end we show that for small $\mu$ the map $F_\mu$ is a
contraction in some ball centered at $- \frac{a_1(y_0)}{a_2(y_0)}$ if
$a_1(y_0)\neq 0$ and in $B_r(-\frac{a_3(y_0)}{a_2(y_0)}\mu)$, if $a_1(y_0)=0$, where 
$$0<r\le r_0=r_0(\mu):=\frac 12\mu \frac{|a_3(y_0)|}{|a_2(y_0)|}.$$ 
To prove that $F_\mu$ is a contraction we may proceed as in Lemma \ref{l:implicit:curve}.
We only need the derivative of $\gamma$, which is given in Lemma \ref{l:d_gamma_d_eps}. 
\end{proof}


\section{A priori estimates}
\label{s:a priori_estimates}
We combine the results of Sections \ref{s:blow_up}-\ref{s:solv_alpha_equal_null} to prove
the $C^2$-a priori estimates announced in the introduction.
\begin{theorem}
\label{t:a priori_ext}
Suppose there is $A_0>2$ such that $k\in C^5(S^3)$ satisfies 
\begin{align*}
D^2 k_\theta(0) \text{ is invertible, if }  
\theta \in \mathcal{A}:= \{\theta \in S^3 \where \nabla k(\theta)=0 \text{ and }\laplace k(\theta)=0\}, 
\end{align*}
\begin{align*}
(A_0)^{-1} \le 1+(1+A_0^{-1})k(\theta) \le A_0,
\end{align*}
$\|k\|_{C^5(S^3)}\le A_0$, and
\begin{align*}
A_0 \ge \supm\{\big\|\big(D^2k_\theta(0)\big)^{-1}\big\| \where \theta \in \mathcal{A}\}.\\
\end{align*}
Thus, $\mathcal{A}$ is discrete and there is $r=r(A_0)>0$ such that
\begin{align*}
\nabla k(\theta) \neq 0 \text{ for all } \theta \in \cupl_{\theta_0 \in \mathcal{A}}
\overline{B_r(\theta_0)} \setminus \{\theta_0\}.  
\end{align*}
Additionally, assume there is $A_1>0$ such that
\begin{align*}
A_1\ge \supm\{|\laplace k_\theta(0)|^{-1} \where |\nabla k(\theta)|\le A_1^{-1} \text{ and } 
\theta \in S^3 \setminus \cupl_{\theta_0 \in \mathcal{A}}B_r(\theta_0) \},\\
A_1\ge \supm\{|a_0(\theta)|^{-1} \where \theta \in \mathcal{A} \text{ and } a_0(\theta)\neq 0\}.  
\end{align*}
Denote by $M$ the finite set
\begin{align*}
M:= \big\{\theta \in S^3 \where \theta \in \mathcal{A}, \, a_0(\theta)= 0, \,\text{and } a_2(\theta)\neq 0\big\}. 
\end{align*}
Then for every $\delta>0$ exits $C=C(\mathcal{A},A_0,A_1,\delta)$ such that for all 
\begin{align*}
t \in (0,1]\setminus \cup_{\theta \in M} B_\delta(-a_1(\theta)/a_2(\theta))  
\end{align*}
and solutions $\phi_t$ of \eqref{eq:eq1} we have
\begin{align*}
C^{-1} \le \phi_t(x) \le C  \ \text{ and } \|\phi_t(x)\|_{C^{2,\a}(S^3)}\le C.
\end{align*}
\end{theorem}
\begin{proof}
Set $I_{\delta,k} := (0,1]\setminus \cup_{y \in M} B_\delta(-a_1(y)/a_2(y))$. 
To obtain a contradiction, we assume that there are sequences $(k_i) \in C^5(S^3)$, satisfying the assumptions
of the theorem with $(\mathcal{A},A_0,A_1,\delta)$ fixed, and  
$(t_i,\phi_{t_i}) \in I_{\delta,k_i}\times C^2(S^3)$ of solutions
to \eqref{eq:eq1} with $k=k_i$ such that $\|\phi_{t_i}\|_\infty \to \infty$ as $i \to \infty$.
Passing to a subsequence we may assume $t_i \to t_0$ as $i \to \infty$.
By Corollary \ref{c:blow_up} there are $\theta \in S^3$, $\mu_i \to 0$ and $y_i \to 0$
such that $u_{t_i}$ defined by \eqref{eq:116} in stereographic coordinates 
$\stereo_{\theta}(\cdot)$ solves \eqref{eq:eq2} and satisfies
\begin{align*}
{\tilde w}_{t_i}:= u_{t_i}-(1+t_i (k_i)_\theta (y_i))^{-\frac14} z_{\mu_i,y_i} \text{ is orthogonal to }
T_{z_{\mu_i,y_i}}Z,\\
\|{\tilde w}_{t_i}\|_{\Did} = o(1).
\end{align*}
Using the notation of Lemma \ref{p:implicit} we have with $k=k_i$
\begin{align*}
0= f_{t_i}'(u_{t_i})= f_{t_i}'(z_{\mu_i,y_i}+w_0(t_i,\mu_i,y_i)+ {\tilde w}_{t_i}).  
\end{align*}
Consequently, for large $i$, due to the uniqueness of $\vec{\a}$ and $w$ in Lemma \ref{p:implicit},
\begin{align*}
u_{t_i}=z_{\mu_i,y_i}+ w(t_i,\mu_i,y_i,k_i) \text{ and } \vec{\a}(t_i,\mu_i,y_i,k_i)=0,  
\end{align*}
where we added the additional parameter $k_i$ to express the dependence of $\vec{\a}$ and $w$ on $k_i$.
From the expansion of $\vec{\a}$ in \eqref{eq:14} we see
\begin{align*}
\lim_{i \to \infty} \nabla (k_i)_\theta(y_i)=0 
\text{ and } \lim_{i \to \infty} \laplace (k_i)_\theta(y_i)=0.  
\end{align*}
As $(\mathcal{A},A_0,A_1,\delta)$ is fixed, the point $\theta$ is in $\mathcal{A}$, hence $\theta$ is a nondegenerated
critical point of each $k_i$. We may apply Lemma \ref{l:implicit:curve} with $k=k_i$
and get for large $i$
\begin{align*}
y_i = \beta(t_i,\mu_i,k_i),  
\end{align*}
where again the additional parameter $k_i$ denotes the dependence on $k_i$.
From Lemma \ref{l:expansion_tilde_alpha_0} we now get
\begin{align}
\label{eq:121}
0 &= \frac{1}{t_i \mu_i^3} \Big(\vec{\a}\big(t_i,\mu_i,\beta(t_i,\mu_i,k_i),k_i\big)\Big) \notag \\
&= -(1+t_i k_i(\theta))^{-\frac54} \frac{3^\frac34 4}{\pi \sqrt{5}} a_0(\theta,k_i) \notag \\
&\quad +\mu_i (1+t_i k_i(\theta))^{-\frac94} \frac{\pi 3^\frac34 \sqrt{5}}{30}
\Big(a_1(\theta,k_i)+t_i a_2(\theta,k_i)\Big) \notag \\
&\quad + O(\mu_i^{1+\frac14}).
\end{align}
Consequently, as $(\mathcal{A},A_0,A_1,\delta)$ is fixed, 
$a_0(\theta,k_i)=0$ for large $i$.\\ 
We observe that $|a_2(\theta,k_i)| \ge A_0^{-4}$ for all $i$ large enough, if not
then we get, up to a subsequence, 
\begin{align*}
|a_1(\theta,k_i)| &\ge  |k_i(\theta)|^{-1} \bigg(\frac{15}{8\pi}
\Big|\intg_{\rand B_1(0)}\big|D^2(k_i)_{\theta}(0)(x)^2\big|^2\Big|-A_0^{-4}\bigg)  \\
&\ge {\rm const} A_0^{-3}(1-A_0^{-2}),
\end{align*}
which yields a positive lower bound on $|a_1(\theta,k_i)+t_i a_2(\theta,k_i)|$
contradicting the expansion in \eqref{eq:121} for $i$ large. 
Hence from \eqref{eq:121} we infer
\begin{align*}
\Big|t_i + \frac{a_1(\theta,k_i)}{a_2(\theta,k_i)}\Big| \le |a_2(\theta,k_i)|^{-1}O(\mu_i^{\frac14})
\le O(\mu_i^{\frac14}),  
\end{align*}
which is impossible for $\delta >0$. This shows that all solutions $\phi_{t}$ of \eqref{eq:eq1} 
with $t\in I_{\delta}$ are uniformly bounded. From Harnack's inequality and standard elliptic
estimates they are uniformly bounded below by a positive constant and uniformly bounded in $C^{2,\a}(S^3)$,
which ends the proof.
\end{proof}
From the proof of Theorem~\ref{t:a priori_ext} it is clear that $k$ need only to be in $C^4(S^3)$, but then the
constant $C$ will also depend on the modulus of continuity of $D^4 k$. 

\begin{proof}[Proof of Theorems \ref{t:a priori_existence_plus} and \ref{t:a priori_existence_null}]
If $\theta \in M^*_+ \cup M^*_0$ we may apply 
Lemmas \ref{l:implicit:curve} and \ref{l:implicit:t(mu)} with $k=k_\theta$ and $y_0=0$.
If we set $y(\mu):= \beta\big(\tilde{t}(\mu),\mu\big)$ then we have 
$\vec{\a}\Big(\tilde{t}(\mu),\mu,y(\mu)\Big)=0$
for all $0 < \mu < \min(\delta_1,\delta_2)$ and
$
y(\mu)= O(\mu^2).
$
From Lemma \ref{p:implicit} we get that
$$\psi(\mu):= z_{\mu,y(\mu)}+w(\tilde{t}(\mu),\mu,y(\mu))$$ 
is a solution of \eqref{eq:eq2} with $t=\tilde{t}(\mu)$.
As  $\nabla k_\theta(0)=0$ and $y(\mu)= O(\mu^2)$ we may use \eqref{eq:105}
to obtain in $\Did$
\begin{align*}
\psi(\mu) 
&= (1+\tilde{t}(\mu)k_\theta(0))^{-\frac14} z_{\mu,0} +O(\mu^2).  
\end{align*}
To show that $\psi(\mu)$ is positive for small $\mu$, we note that from Sobolev's inequality
$\psi(\mu)^{-} \to 0$ in $L^6$ as $\mu \to 0$, where $\psi(\mu)^{-}:= \min(\psi(\mu),0)$.
Testing $f_t'(\psi(\mu))$ with $\psi(\mu)^{-}$ and using Sobolev's inequality we get for some $c(k)>0$
\begin{align*}
\int |\nabla\psi(\mu)^{-}|^2 = \int (1+\tilde{t}(\mu)k_\theta(x)) (\psi(\mu)^{-})^6
\le c(k) \Big(\int |\nabla\psi(\mu)^{-}|^2\Big)^3.
\end{align*}
If $\psi(\mu)^{-} \neq 0$ for small $\mu$ we obtain the contradiction
\begin{align*}
c(k)^{-\frac12} \le \int |\nabla\psi(\mu)^{-}|^2 = \int (1+\tilde{t}(\mu)k_\theta(x)) (\psi(\mu)^{-})^6
\xpfeil{\mu \to 0}0.  
\end{align*}
The $C^0$-estimate then follows from elliptic regularity (see \cite{BrezisKato79}). 
Setting 
$$\phi^\theta(\mu,\cdot):= (\mathcal{R}_\theta)^{-1}(\psi(\mu)) \text{ and } t^\theta(\mu)=\tilde{t}(\mu)$$  
yields the existence of the desired curve of solutions.\\
To prove uniqueness of the curves suppose $(t_i,\phi_i)$ blows up at $\theta \in S^3$.
If $t_i \in (\delta,1+\delta)$ then, as in the proof of Theorem \ref{t:a priori_ext}, we get $\theta \in M^*_+$.
Under the assumptions of Theorem \ref{t:a priori_existence_null} we already know that $\theta \in M^*_0$.
If all but finitely many $(t_i,\phi_i)$ lie on the curve corresponding 
to $\theta \in M^*_+ \cup M^*_0$, we are done. Hence we may assume,
going to a subsequence if necessary, that none of the $(t_i,\phi_i)$ lie on the curve. This
is impossible since by Corollary \ref{c:blow_up} and Lemma \ref{p:implicit}
there are $\mu_i,y_i$ converging to zero such that for $i$ large
\begin{align*}
\mathcal{R}_\theta(\phi_i) = z_{\mu_i,y_i}+w(t_i,\mu_i,y_i) \text{ and } \vec{\a}(t_i,\mu_i,y_i)=0,
\end{align*}
and thus applying the uniqueness part in Lemmas \ref{l:implicit:curve} and \ref{l:implicit:t(mu)}
we see that $y_i= \beta(t_i,\mu_i)$ and $t_i= \tilde{t}(\mu_i)$ and the points $(t_i,\phi_i)$
have to lie on the curve.  
\end{proof}

\appendix

\section{Formulas and integrals}
\label{s:formulas-integrals}
The {\em Jacobi} polynomial ${\mathcal P}_j^{(\sigma,\sigma)}$ is defined by
\begin{align}
\label{eq:A1}
{\mathcal P}_j^{(\sigma,\sigma)}(x) 
:= \frac{(-1)^j}{2^j j!}(1-x^2)^{-\sigma} \frac{d^j}{dx^j}\big((1-x^2)^{\sigma+j}\big)  
\end{align}
To compute integrals containing Jacobi-polynomials we will use
\begin{align}
\label{eq:A45}
\int_{-1}^{1} (1-\xi^2)^{\sigma} {\mathcal P}_{j}^{(\sigma,\sigma)}(\xi) {\mathcal P}_{i}^{(\sigma,\sigma)}(\xi) 
= \frac{2^{2\sigma+1} \Gamma(j+\sigma+1)^2}{(2j+2\sigma+1)\,j!\, \Gamma(j+2\sigma+1)} \delta_{i,j},
\end{align}
and the recurrence relation for $j \in \nz_0$
\begin{align}
\label{eq:A46}
\begin{split}
\xi {\mathcal P}_{j}^{(\sigma,\sigma)}(\xi) = 
\frac{(j+1)(j+2\sigma+1)}{(2j+2\sigma+1)(j+\sigma+1)} {\mathcal P}_{j+1}^{(\sigma,\sigma)}(\xi) +
\frac{j+\sigma}{2j+2\sigma+1} {\mathcal P}_{j-1}^{(\sigma,\sigma)}(\xi),\\
{\mathcal P}_{-1}^{(\sigma,\sigma)}(\xi):= 0, \; {\mathcal P}_{0}^{(\sigma,\sigma)}(\xi) := 1.    
\end{split}
\end{align}
For a detailed account on {\em Jacobi} polynomial we refer to \cite{szego75}.
To evaluate integrals of the form 
\begin{align*}
\int_0^\infty r^a(1+r^2)^{-b} {\mathcal P}_{j}^{(\sigma,\sigma)}\Big(1-\frac{2}{1+r^2}\Big) \di r,
\end{align*}
we use the following change of coordinates
\begin{align}
\label{eq:A44}
\xi= \frac{r-r^{-1}}{r+r^{-1}} = 1- \frac{2}{1+r^2},  
\end{align}
which gives
\begin{align*}
r = \Big(\frac{1+\xi}{1-\xi}\Big)^{1/2} \text{ and } dr = \Big(\frac{1+\xi}{1-\xi}\Big)^{-1/2}(1-\xi)^{-2} d\xi,  
\end{align*}
and leads to 
\begin{align*}
2^{-b}\int_{-1}^1 (1+\xi)^{\frac{a-1}{2}} (1-\xi)^{b-2-\frac{a-1}{2}} {\mathcal P}_{j}^{(\sigma,\sigma)}(\xi) \di \xi.  
\end{align*}
Moreover, we note that for any $a,b>-1$ 
\begin{align}
\label{eq:A40}
\int_{-1}^{1}(1+\xi)^a (1-\xi)^b d\xi= 2^{a+b+1} \frac{\Gamma(a+1) \Gamma(b+1)}{\Gamma(a+b+2)}.     
\end{align}
This gives for $a>-1$ and $2b-a>1$
\begin{align}
\label{eq:A34}
\int_0^\infty r^a (1+r^2)^{-b} = \frac{\Gamma(1+(a-1)/2) \Gamma(b-1-(a-1)/2)}{2 \Gamma(b)}. 
\end{align}
To compute integrals over $\rz^N$ we use polar coordinates. To compute the resulting integrals
over $\rand B_1(0)$ we use the following elementary results:\\
For $\vec{\beta}\in \nz_0^N$ we have
\begin{align}
\label{eq:A8}
\int_{\rand B_1(0)} \prod_{i=1}^N x_i^{2\beta_i} = 
\frac{2\prod_{i=1}^N \Gamma(\beta_i+\frac12)}{\Gamma\big(\frac{N}2+\sum_{i=1}^N \beta_i\big)}.  
\end{align}
Let $m \ge 2$ and $P_m$ be a homogeneous polynomial of order $m$ in $x \in \rz^N$. Then
\begin{align}
\label{eq:A33}
\intg_{\rand B_1(0)} P_m(x) &= \frac{1}{2N+(N+m)(m-2)}\intg_{\rand B_1(0)} (\laplace P_m)(x),\\
\label{eq:A70}
\intg_{\rand B_1(0)} P_m(x) x_i &= 
\frac{1}{2(N+2)+(N+m+1)(m-3)}\intg_{\rand B_1(0)} (\laplace P_m)(x) x_i.
\end{align}
As $\frac{\rand}{\rand x_i} D^\ell_y k(y)(x)^\ell = \ell \frac{\rand}{\rand y_i} D^{\ell-1}_y k(y)(x)^{\ell-1}$
we see
\begin{align}
\label{eq:A35}
\intg_{\rand B_1(0)} \frac{D^{2\ell}k(y)(x)^{2\ell}}{(2\ell)!} &=  
\frac{2 \pi^{N/2}\laplace^\ell k(y)}{\Gamma(N/2) \prod_{m=1}^\ell 2m(N+2m-2)},
\end{align}
\begin{align}
\label{eq:A36}
\intg_{\rand B_1(0)} \frac{D^{2\ell+1}k(y)(x)^{2\ell+1} x_i}{(2\ell+1)!}&=  
\frac{2 \pi^{N/2} \frac{\rand}{\rand y_i}\laplace^\ell k(y)}{N\Gamma(N/2) \prod_{m=1}^\ell 2m(N+2m)}.
\end{align}

\begin{lemma}
\label{l:A3}
Suppose $j \in \nz_0$, then
\begin{align}
\label{eq:A43}
\int_0^\infty r^6 (1+r^2)^{-5}\Big(1-\frac{2}{1+r^2}\Big) 
{\mathcal P}_j^{(\frac52,\frac52)}\Big(1-\frac{2}{1+r^2}\Big)  \notag \\
= \frac{\Gamma(\frac32)\Gamma(\frac72+j) (j^2+6j+2)}{2 \Gamma(6+j)}.
\end{align}
\end{lemma}
\begin{proof}
We use the change of coordinates in \eqref{eq:A44} and obtain
\begin{align*}
\int_0^\infty &r^4 (1+r^2)^{-5}\Big(1-\frac{2}{1+r^2}\Big) 
{\mathcal P}_j^{(\frac52,\frac52)}\Big(1-\frac{2}{1+r^2}\Big)\\
&= 2^{-5} \int_{-1}^1 (1-\xi^2)^{\frac52} \frac{\xi}{(1-\xi)^2}{\mathcal P}_j^{(\frac52,\frac52)}(\xi)\\
&= \frac{(-1)^j}{2^{j+5}j!} \int_{-1}^1 \Big(\frac{1}{(1-\xi)^2}- \frac{1}{1-\xi}\Big) 
\frac{d^j}{d \xi^j}\big((1-\xi^2)^{\frac52+j}\big)\\
&= \frac{(-1)^j}{2^{j+5}j!} \int_{-1}^1 
\Big(\frac{(-1)^j(j+1)!}{(1-\xi)^{2+j}}- \frac{(-1)^j j!}{(1-\xi)^{1+j}}\Big) 
(1-\xi^2)^{\frac52+j}.
\end{align*}
Now, the claim follows from \eqref{eq:A40}. 
\end{proof}
\begin{lemma}
\label{l:A4}
Let $\beta_j :=  \frac{(3+j)j!}{\Gamma(j+\frac{7}2)}$ for $j \in \nz$. Then  
\begin{align}
\label{eq:A19}
\int_0^\infty r^6 (1+r^2)^{-7} \Big(1-\frac{2}{1+r^2}\Big) 
\bigg(\sum_{j=0}^\infty \beta_j {\mathcal P}_j^{(\frac52,\frac52)}\Big(1-\frac{2}{1+r^2}\Big)\bigg)^2
= \frac{1}{288}.
\end{align}
\end{lemma}
\begin{proof}
We use the change of variable given in \eqref{eq:A44} and the recurrence formula \eqref{eq:A46},
applied to $\xi {\mathcal P}_j^{(\frac52,\frac52)}(\xi)$, and get
\begin{align*}
\int_0^\infty r^6 &(1+r^2)^{-7} \Big(1-\frac{2}{1+r^2}\Big) 
\bigg(\sum_{j=0}^\infty \beta_j {\mathcal P}_j^{(\frac52,\frac52)}\Big(1-\frac{2}{1+r^2}\Big)\bigg)^2  \\
&= 2^{-7} \int_{-1}^1 (1-\xi^2)^{\frac52} \xi 
\bigg(\sum_{j=0}^\infty \beta_j {\mathcal P}_j^{(\frac52,\frac52)}(\xi)\bigg)^2\\
&= 2^{-7} \sum_{d=0}^{\infty} \sum_{j=0}^d \beta_j \beta_{d-j} 
\int_{-1}^1 (1-\xi^2)^{\frac52} \xi {\mathcal P}_j^{(\frac52,\frac52)}(\xi)  
{\mathcal P}_{d-j}^{(\frac52,\frac52)}(\xi)\\
&=  \sum_{l=0}^{\infty} \beta_l \beta_{l+1} 
\frac{\Gamma(l+\frac72)^2 (l+\frac72)}{\Gamma(l+6)(2l+8)(2l+6) l!}
= \frac{1}{4}\sum_{l=0}^{\infty} \frac{(l+1)!}{(l+5)!} = \frac{1}{288}.
\end{align*}  
\end{proof}

\bibliographystyle{adinat}
\bibliography{prescribed}

\end{document}